\documentclass{amsart}

\usepackage[T1]{fontenc}
\usepackage[utf8]{inputenc}
\usepackage{float}
\usepackage{setspace}
\usepackage{amsmath,amsfonts,amssymb,amsthm,mathrsfs}
\usepackage{graphicx,psfrag,epsf}
\usepackage{enumerate}
\usepackage{url} % not crucial - just used below for the URL 

\DeclareMathOperator{\dH}{d}
\DeclareMathOperator{\pen}{pen}
\DeclareMathOperator{\KL}{KL}
\DeclareMathOperator{\crit}{crit}

 \usepackage{import} 

\usepackage[english]{babel}

\newtheorem{defi}{Definition}[section]

\newtheorem{theorem}[defi]{Theorem}
\newtheorem{corollary}{Corollary}
\newtheorem{lemma}[defi]{Lemma}

\newtheorem{assum}{Assumption}
\newtheorem{proposition}{Proposition}
\newenvironment{demo}{\noindent{\bf Proof.}}{\qed}
\newcommand{\y}{\mathbf{y}}

\newcommand{\SSigma}{\boldsymbol{\Sigma}}

\title[ Block-diagonal covariance selection for high-dimensional GGM ]{Block-diagonal covariance selection for high-dimensional Gaussian graphical models}% At most 5 thanks
\author{Emilie Devijver}\address{Department of Mathematics and Leuven Statistics Research Center (LStat), KU Leuven, Leuven, Belgium}
\author{M\'elina Gallopin}\address{Laboratoire MAP5, Universit\'e Paris Descartes and CNRS, Sorbonne Paris Cit\'e, Laboratoire de Math\'ematiques, UMR 8628, B\^atiment 425, Universit\'e Paris-Sud, F-91405, Orsay, France, INRA, UMR 1313 G\'en\'etique animale et biologie int\'egrative, 78352 Jouy-en-Josas, France}

\date{\today}

\begin{document}

 \begin{abstract}
Gaussian graphical models are widely utilized to infer and visualize networks of dependencies between continuous variables. However, inferring the graph is difficult when the sample size is small compared to the number of variables. To reduce the number of parameters to estimate in the model, we propose a non-asymptotic model selection procedure supported by strong theoretical guarantees based on an oracle type inequality and a minimax lower bound. The covariance matrix of the model is approximated by a block-diagonal matrix. The structure of this matrix is detected by thresholding the sample covariance matrix, where the threshold is selected using the slope heuristic. Based on the block-diagonal structure of the covariance matrix, the estimation problem is divided into several independent problems: subsequently, the network of dependencies between variables is inferred using the graphical lasso algorithm in each block. The performance of the procedure is illustrated on simulated data. An application to a real gene expression dataset with a limited sample size is also presented: the dimension reduction allows attention to be objectively focused on interactions among smaller subsets of genes, leading to a more parsimonious and interpretable modular network.

 \end{abstract}

\maketitle
\tableofcontents

 \section{Introduction}
 
%% INTRODUCTION
Graphical models \cite{Whittaker1990} have become a popular tool for representing conditional dependencies among variables using a graph. For Gaussian graphical models (GGMs), the edges of the corresponding graph are the non-zero coefficients of the inverse covariance matrix. To estimate this matrix in high-dimensional contexts, methods based on an $\ell_1$ penalized log-likelihood have been proposed \cite{MeinshausenBuhlmann, Yuan2007, Banerjee2008}. A popular method is the graphical lasso algorithm introduced by \cite{Friedman2007}. Gaussian graphical models have many potential applications, such as the reconstruction of regulatory networks from real omics data \cite{Krumsiek2011,Akbani2014pan}. However, these methods often perform poorly in so-called \textit{ultra high-dimensional contexts} \cite{Giraud2008, Verzelen2012}, when the number of observations is much smaller than the number of variables. A small sample size is a common situation in various applications, such as in systems biology where the cost of novel sequencing technologies may limit the number of available observations \cite{Frazee2011}. In practice, the network reconstruction problem is facilitated by restricting the analysis to a subset of variables, based on external knowledge and prior studies of the data \cite{2009_EJS_Chiquet,Yin2011}. When no external knowledge is available, only the most variable features are typically kept in the analysis \cite{Guo2011,Allen2013}. Choosing the appropriate subset of variables to focus on is a key step in reducing the model dimension and the number of parameters to estimate, but no procedure is clearly established to perform this selection when the sample size is really low.

% In this paper, we propose a non-asymptotic data-driven method to detect groups of variables with relevant interactions.

In the context of graphical lasso estimation, \cite{Witten2011GM} and \cite{Mazumder2012} have noticed a particular 
property: the block-diagonal structure of the graphical lasso solution is totally determined by the block-diagonal structure of the thresholded empirical covariance matrix. The graphical lasso estimation for a given level of regularization $\lambda$ can be decomposed into two steps. First, the connected components of the graphical model are detected, based on the absolute value of the sample covariance matrix thresholded at level $\lambda$. Second, the graphical lasso problem is solved in each connected component independently using the same regularization parameter $\lambda$ for each subproblem. This decomposition is of great interest to reduce the number of parameters to estimate for a fixed level of regularization. It has been exploited for large-scale problems \cite{Zhao2012} and for joint graphical lasso estimations \cite{Danaher2014}. \cite{Hsieh2014} have improved the computational cost of the two-step procedure by using a quadratic approximation, and have proved the superlinear convergence of their algorithm. \cite{Tan2015} have noticed that the first step, {\it i.e.} the detection of connected components by thresholding, is equivalent to performing a single linkage clustering on the variables. They have proposed the \textit{cluster graphical lasso}, using an alternative to single linkage clustering in the two-step procedure, such as average linkage clustering. The selection of the cutoff applied to hierarchical clustering in the first step of the \textit{cluster graphical lasso} algorithm is performed independently from the  selection of the regularization parameters in the second step of the algorithm. Their results suggest that the detection of the block diagonal structure of the covariance matrix prior to network inference in each cluster can improve network inference. Other authors have recently proposed procedures to detect the block-diagonal structure of a covariance matrix. \cite{Pavlenko2012} provided a method to detect this structure for high-dimensional supervised classification that is supported by asymptotic guarantees. \cite{Hyodo2015} proposed tests to perform this detection and derived consistency for their method when the number of variables and the sample size tend to infinity.

 In this paper, we adapt the two-step procedure proposed by \cite{Witten2011GM} and \cite{Mazumder2012} to infer networks of conditional dependencies between variables. To improve network inference performance, we decouple the two steps of the procedure, using different parameters for thresholding the empirical covariance matrix and for estimation within each connected component of the network, as  proposed by \cite{Tan2015}. The main contribution of our work stands in the use of a non-asymptotic criterion to perform the first step of the two-step procedure, {\it i.e.} the detection of the block-diagonal structure of a covariance matrix. In our procedure, we recast the detection problem into a model selection problem and aim at choosing the best model among a collection of multivariate distributions with block-diagonal covariance matrices. We obtained non-asymptotic theoretical guaranties, based on the control of the risk of the model selection procedure. These results provide a penalty to select a model. This penalty is known up to multiplicative constants depending on assumptions hard to explicitly satisfy in practice. To calibrate these constants in practice, we use the slope heuristic, originally proposed by \cite{BirgeMassart2001} and detailed in \cite{Baudry2012}. Unlike other methods to detect the appropriate block-diagonal covariance matrix \cite{Pavlenko2012,Tan2015,Hyodo2015}, our procedure is non-asymptotic and offers strong theoretical guarantees when the number of observations is limited, which is of great interest for many real applications. More precisely, we prove that our estimator is adaptive minimax to the structure of the covariance matrix.
 
The paper is organized as follows. In Section~\ref{OurMethod}, after providing basic notations and definitions, the non-asymptotic method to detect the block-diagonal structure of the GGM is 
presented, as well as the complete framework to infer network where the number of observations is limited. Section~\ref{TheoreticalResults} details theoretical results supporting our model 
selection criterion. In particular, an oracle type inequality upper bounds the risk between the true model and the model 
selected among the collection of models, and a minimax lower bound guarantees that the non-asymptotic procedure has an optimal rate of convergence.
Section~\ref{Simul} investigates the numerical performance of our method in a simulation study. Section~\ref{RealData} illustrates our procedure on a real gene expression RNA-seq dataset with 
a limited sample size. 
{All proofs are provided in Appendix 1.}
 %%%%%%%%%%%%%%%%%%%%%%%%%%%%%%%%%%%%%%%%%%%%%
%% METHOD
%%%%%%%%%%%%%%%%%%%%%%%%%%%%%%%%%%%%%%%%%%%%%

 \section{A method to detect block-diagonal covariance structure}
\label{OurMethod}

Let $\y = (\y_1,\ldots,\y_n)$ be a sample in $\mathbb{R}^p$ from a multivariate normal distribution with density $\phi_p(0,\Sigma)$ where $\Sigma_{j,j}=1$ for all 
$j \in \{1,\ldots,p\}$. Let $S$ be the empirical covariance matrix associated with the sample $\y$. Our goal is to infer the graph of conditional dependencies between variables, encoded by 
the precision matrix $\Theta=\Sigma^{-1}$. Since the matrices $\Sigma$ and $\Theta$ have the same block-diagonal structure, we first seek to detect the optimal block-diagonal structure of the 
covariance matrix ${\Sigma}$, {\it i.e.} the optimal partition of variables into blocks. {We index the variables from $1$ to $p$. 
We note $\mathbf{B}=\{\mathbf{B}_1; \ldots; \mathbf{B}_K\}$ the partition of variables into $K$ blocks where $\mathbf{B}_k$ is the subset of variables in block $k$, and $p_k$ is the number 
of variables in block $k$. The partition describes the block-diagonal structure of the matrix: off the block, all coefficients of the matrix are zeros. Note that we authorize the permutation 
inside blocks: e.g. the partition of $5$ variables into blocks $\{(1,3,4); (2,5)\}$ is equivalent to the partition $\{(3,1,4); (5,2)\}$. We also authorize reordering of the blocks: e.g. 
$\{(1,3,4); (2,5)\}$ is equivalent to $\{(2,5);(1,3,4)\}$. We consider the following set of multivariate normal densities with block-diagonal covariance matrices:} 
{
\begin{equation}
F_\mathbf{B} = \left\{ f_\mathbf{B} =\phi _p(0,{{\Sigma}}_\mathbf{B}) \text{ with } \Sigma_\mathbf{B} \in \mathbb{S}_p^{++}(\mathbb{R}) \left| 
\begin{array}{l} \lambda_m \leq {\min}(\text{sp}({\Sigma}_\mathbf{B})) \leq {\max}(\text{sp}({\Sigma}_\mathbf{B})) \leq \lambda_M, \\
% \lambda_m \leq \Lambda_{\min} (\Sigma_\mathbf{B}) \leq \Lambda_{\max}(\Sigma_\mathbf{B}) \leq \lambda_M,\\
 \Sigma_\mathbf{B} = P_\sigma 
\begin{pmatrix}
 {\Sigma}_1 & 0 & 0\\
 0 & \ddots & 0 \\
 0 & 0 & {\Sigma}_K
\end{pmatrix}P_\sigma^{-1}, \\
 {\Sigma}_k \in \mathbb{S}_{p_k}^{++}(\mathbb{R}) \text{ for } k \in \{1,\ldots,K\} \\
\end{array}\right. \right\}, \\
%S_\mathbf{B} &= \left\{ \Sigma_\mathbf{B} \in \mathbb{S}_p^{++}(\mathbb{R})\left| 
%\begin{array}{l} e_m \leq {\min}({\Sigma}_\mathbf{B}) \leq {\max}({\Sigma}_\mathbf{B}) \leq e_M, \\
% \lambda_m \leq \Lambda_{\min} (\Sigma_\mathbf{B}) \leq \Lambda_{\max}(\Sigma_\mathbf{B}) \leq \lambda_M,\\
% \Sigma_\mathbf{B} = P_\sigma 
%\resizebox{2cm}{!}{$\begin{pmatrix}
% {\Sigma}_1 & 0 & 0\\
% 0 & \ddots & 0 \\
% 0 & 0 & {\Sigma}_K
%\end{pmatrix}$}P_\sigma^{-1}, \\
% {\Sigma}_k \in \mathbb{S}_{p_k}^{++}(\mathbb{R}) \text{ for } k \in \{1,\ldots,K\}, \\
%\end{array}
%\right\},
%%\right. \right. \nonumber \\
%% &\left. \phantom{\begin{pmatrix}
%% {\Sigma}^1 & 0 & 0\\
%% 0 & \ddots & 0 \\
%% 0 & 0 & {\Sigma}^K
%%\end{pmatrix} } {\Sigma}_k \in \mathbb{S}_{p_k}^{++}(\mathbb{R}) \text{ for } k \in \{1,\ldots,K\} \right\}. 
\label{myFB}
\end{equation}
where $\mathbb{S}_p^{++}(\mathbb{R})$ is the set of positive semidefinite matrices of size $p$, $\lambda_m$ and $\lambda_M$ are real numbers, ${\min}(\text{sp}({\Sigma}_\mathbf{B})), {\max}(\text{sp}({\Sigma}_\mathbf{B}))$ 
are the smallest and highest eigenvalues of $\Sigma_\mathbf{B}$ and $P_\sigma$ is a permutation matrix leading to a block-diagonal covariance matrix.
}

{We consider $\mathcal{B}$ the set of all possible partitions of variables. In theory, we would like to consider the collection of models 
 \begin{align}
 \label{modelCollectionFmathcal}
\mathcal{F} = (F_\mathbf{B})_{\mathbf{B} \in \mathcal{B} }. 
 \end{align}
 However, the set $\mathcal{B}$ is large: its size is the Bell number. An exhaustive exploration of the set $\mathcal{B}$ is then not possible even for a moderate number of variables $p$.
% $$S(p,k) =\frac{1}{k!} \sum_{j=0}^k (-1)^{k-j} \dbinom{j}{k}j^p.$$
 We restrict our attention to the sub-collection:
 \begin{align}
 \label{Blambda}
 \mathcal{B}_{\Lambda}=(\mathbf{B}_\lambda)_{\lambda \in \Lambda}
 \end{align}
 of $\mathcal{B}$ where $\mathbf{B}_\lambda$ is 
 the partition of variables corresponding to the block-diagonal structure of the adjacency matrix $E_{\lambda}=[\mathbf{1}_{\{ \mid S_{j,j^\prime}\mid > \lambda\} }]_{ 1 \leq j \leq p \atop 1 \leq j^\prime \leq p }$,
 based on the thresholded absolute value of the sample covariance matrix $S$.
Recall that \cite{Mazumder2012} have proved that {the class of block-diagonal structures $\mathcal{B}_{\Lambda}$ detected by thresholding of the sample covariance is the same class of block-diagonal structures
detected by the graphical lasso algorithm when the regularization parameter varies, which supports the fact that we restrict our attention to this specific sub-collection.}
 Note that the data is scaled if needed so that the set of thresholds $\Lambda \subset [0,1]$ covers all possible partitions derived 
 from $E_{\lambda}$.}
 
{Once we have constructed the collection of models $\mathcal{F}_\Lambda = (F_\mathbf{B})_{\mathbf{B} \in \mathcal{B}_\Lambda}$, we select a model among this collection using the following model selection criterion:
\begin{align}
\label{penkappaDm}
\hat{\mathbf{B}} &= \underset{\mathbf{B} \in \mathcal{B}_{\Lambda}}{\operatorname{argmin}} \left\{ -\frac{1}{n} \sum_{i=1}^n \log(\hat{f}_\mathbf{B}(\y_i)) + \mathrm{pen}(\mathbf{B}) \right\}, \nonumber
%\mathrm{pen}(\mathbf{B})&=\kappa D_\mathbf{B}, \\
%\mathrm{pen}(\mathbf{B})&=\kappa D_\mathbf{B} 
%\kappa_1 D_\mathbf{B}/n + \kappa_2 D_\mathbf{B}/n log(p^ 2/(D_\mathbf{B}/n))
%\mathrm{pen}(\mathbf{B})&=\kappa \frac{D_\mathbf{B}}{n} \left[2c^2 + \log\left( \frac{p^3(p-1)}{D_\mathbf{B}(\frac{D_\mathbf{B}}{n} c^2 \wedge 1)} \right) \right],\\
%\mathrm{pen}(\mathbf{B})&=\kappa_1 \frac{D_\mathbf{B}}{n} + \kappa_2 \frac{D_\mathbf{B}}{n} \log\left( \frac{p^2}{D_\mathbf{B}/n} \right),
\end{align}
}

{where $\mathrm{pen}(\mathbf{B})$ is a penalty term to define and $\hat f_\mathbf{B}=\phi _p(0,\hat{{\Sigma}}_\mathbf{B})$ where $\hat{\Sigma}_\mathbf{B}$ is the maximum 
likelihood estimator of $\Sigma_\mathbf{B}$. 
The matrix $\hat{\Sigma}_\mathbf{B}$ is constructed block by block, using the sample covariance matrix of the dataset restricted to variables in each block. }
 
{The penalty term $\mathrm{pen}(\mathbf{B})$ is based on non-asymptotic model selection properties, as detailed in Section~\ref{TheoreticalResults}: $\mathrm{pen}(\mathbf{B})=\kappa \frac{D_\mathbf{B}}{n} + \tilde{\kappa} \frac{D_\mathbf{B}}{n} \log\left( \frac{p(p-1)}{D_\mathbf{B}} \right)$ where $D_\mathbf{B} = \sum_{k=1}^K {p_k (p_k-1)}/{2}$ is the dimension of the model $F_\mathbf{B}$ and $\kappa, \tilde{\kappa}$ are two constants depending on absolute constants and on the bounds $\lambda_m$ and $\lambda_M$. In practice, we consider a simpler version of the penalty term: 
\begin{equation}
\mathrm{pen}(\mathbf{B})= \kappa \frac{D_\mathbf{B}}{n}
\label{kappaPP}
\end{equation}
}
{where $\kappa$ is a constant depending on absolute constants and on the bounds $\lambda_m$ and $\lambda_M$. Such simplification has already been proposed by \cite{Lebarbier2005}. The extra term 
$\frac{D_\mathbf{B}}{n} \log\left( \frac{p(p-1)}{D_\mathbf{B}} \right)$ is useful to overpenalize the collection of models when it contains many models with the same sizes. The simplification of the 
penalty term \eqref{kappaPP} is reasonable for moderate number of variables. }

{Subsequently, we note that the bounds $\lambda_m$ and $\lambda_M$ are non-tractable. For this reason, we prefer to calibrate the constant $\kappa$ in \eqref{kappaPP} 
from the data. This calibration is based on the slope heuristic, originally proposed and proved in the context of heteroscedastic regression with fixed design \cite{BirgeMassart2007,Baraud2009}, 
and for homoscedastic regression with fixed design \cite{Arlot2010}.} 
{In other contexts, the slope heuristic has been used and have proven to be effective for multiple change point detection \cite{Lebarbier2005}, for variable selection in mixture models \cite{Maugis2011}, for choosing the number of components in Poisson mixture models \cite{Rau2015} or for selecting the number of components in discriminative functional mixture models \cite{Bouveyron2015}. }

\cite{Baudry2012} have provided practical tools to calibrate {the coefficient $\kappa$ in \eqref{kappaPP} based on} the slope heuristic developed by \cite{BirgeMassart2007}. 
We describe these tools in Section \ref{Simul}.
Note that the detection of the optimal $\mathbf{B}$ is easy to implement in practice and does not rely on heavy computation such as cross-validation techniques.

Once we have detected the optimal block-diagonal structure of the GGM, network inference is performed independently in each block using the graphical lasso \cite{Friedman2007}. 

We summarize the method proposed to infer network in high-dimensional context.
\begin{enumerate}
\item[(a)] (Block-diagonal covariance structure detection) Select the modularity structure of the network.
\begin{enumerate}
\item[(1)] Compute the sample covariance matrix $S$.
\item[(2)] Construct the sub-collection of partitions $\mathcal{B}_{\Lambda}=(\mathbf{B}_\lambda)_{\lambda \in \Lambda}$, where $\Lambda$ is a set of thresholds, more precisely the set of values taken by the sample covariance matrix $S$ in absolute value. Each partition corresponds to the block-diagonal structure of the matrix $E_{\lambda}=[\mathbf{1}_{\{ \mid S_{j,j^\prime}\mid > \lambda\} }]_{ 1 \leq j \leq p  \atop 1 \leq j^\prime \leq p }$.
\item[(3)] For each partition $\mathbf{B} \in \mathcal{B}_{\Lambda}$, compute the corresponding maximum log-likelihood of the model. 
\item[(4)] Based on the log-likelihood associated to each partition $\mathbf{B}$ in $\mathcal{B}_{\Lambda}$, calibrate the penalty in equation \eqref{kappaPP} to select the partition $\hat{\mathbf{B}}$ using the slope heuristic.
\end{enumerate}
\item[(b)] (Network inference in each module) For each group of variables in the selected partition $\hat{\mathbf{B}}$, infer the network using the graphical lasso introduced by \cite{Friedman2007}. The choice of the regularization parameter for the graphical lasso algorithm is performed independently in each module. 
\end{enumerate}

%%%%%%%%%%%%%%%%%%%%%%%%%%%%%%%%%%%%%%%%%%%%%
%% THEORETICAL
%%%%%%%%%%%%%%%%%%%%%%%%%%%%%%%%%%%%%%%%%%%%%

 \section{Theoretical results for non-asymptotic model selection}
 \label{TheoreticalResults}

%% INTRO ORACLE INEQUALITY
{
The model selection procedure presented in Section \ref{OurMethod} is justified by theoretical results.
We obtain bounds on the risk between the selected density and the true one, which prove that we select a good block-diagonal structure. More precisely, we aim at selecting, among $\mathcal{B}$, 
the optimal partition $\mathbf{B}^\star$. First, for each model indexed by $\mathbf{B}$, we consider the density $\hat{f}_{\mathbf{B}}=\phi_p(0,\hat{\Sigma}_\mathbf{B})$ where $\hat{\Sigma}_\mathbf{B}$ is the 
maximum likelihood 
estimator of $\Sigma_\mathbf{B}$. Among all $\mathbf{B} \in \mathcal{B}_\Lambda$, we want to select the density $\hat{f}_\mathbf{B}$ which is the closest one to the true distribution $f^\star$. To measure the 
distance between the two densities, we define the risk:
$$ R_\mathbf{B}(f^\star)= \mathbb{E}(d^2(f^\star,\hat{f}^{\mathbf{B}})),$$ where $d$ is a distance between two densities. Ideally, we would like to select the partition $\mathbf{B}$ 
that minimizes the risk $R_\mathbf{B}(f^\star)$: this partition is called the \emph{oracle}. Unfortunately, it is not reachable in practice because the true density $f^\star$ is unknown. 
However, we will prove that we do almost as well as the oracle, i.e. we select a model for which the risk of the procedure is upper bounded by the oracle risk, up to a constant.}

\noindent
% QUELQUES DEFINITIONS AVANT D'ENONCER L'INEGALITE ORACLE 
Before stating the theorem, we recall the definition of the Hellinger distance between two densities $f$ and $g$ defined on $\mathbb{R}^p$,
$\dH_H^2(f,g) = \frac{1}{2} \int_{\mathbb{R}^p} (\sqrt{f(x)}-\sqrt{g(x)})^2 dx,$
and the Kullback-Leibler divergence between two densities $f$ and $g$ defined on $\mathbb{R}^p$,
$\KL (f,g) = \int_{\mathbb{R}^p} \log \left( \frac{f(x)}{g(x)} \right) f(x) dx.$ 
%Besides, the matrix $\Sigma_\mathbf{B}$, defined in \eqref{myFB}, has bounded coefficients $e_m$ and $e_M$. We can prove that it has bounded eigenvalues and denote by $\lambda_m$ and $\lambda_M$ the smallest and the largest eigenvalues of this matrix for practical reasons.We denote $\lambda_m$ and $\lambda_M$ the smallest and the largest eigenvalues of this matrix for practical reasons.

{
\begin{theorem}
\label{inegalite oracle}
 Let $\y=(\y_1,\ldots,\y_n)$ be the observations, arising from a density $f^\star$.
Consider the collection of models $\mathcal{F}$ defined in \eqref{modelCollectionFmathcal}.
We denote by $\hat f_{\mathbf{B}}$ the maximum likelihood estimator for the model $F_\mathbf{B}$. Let $\mathcal{B}_{\Lambda} \subset \mathcal{B}$ as defined in \eqref{Blambda}.
\\
Let $\tau>0$, and for all $\mathbf{B} \in \mathcal{B}$, let $f_\mathbf{B} \in F_\mathbf{B}$ such that: 
\begin{align}
\KL(f^\star,f_\mathbf{B}) &\leq 2 \inf_{f \in F_\mathbf{B}} \KL(f^\star,f) \nonumber ;\\
f_\mathbf{B} &\geq \exp \left({-\tau}\right) f^\star.
\label{Bernstein2}
\end{align}
Then, there exists some absolute constants $\kappa$ and $C_{\text{oracle}}$ such that whenever
 $$\mathrm{pen}(\mathbf{B}) \geq \kappa \frac{D_\mathbf{B}}{n} \left[2c^2 + \log\left( \frac{p^4}{D_\mathbf{B}(\frac{D_\mathbf{B}}{n} c^2 \wedge 1)} \right) \right] $$
 for every $\mathbf{B} \in \mathcal{B}$, with $c = \sqrt{\pi} + \sqrt{\log(3\sqrt{3} \frac{\lambda_M}{\lambda_m})}$, the random variable $\hat{\mathbf{B}} \in \mathcal{B}_\Lambda$ such that
\begin{align}
\label{estimatorB}
\hat{\mathbf{B}} = \underset{\mathbf{B} \in \mathcal{B}_{\Lambda}}{\operatorname{argmin}} \left\{ -\frac{1}{n} \sum_{i=1}^n \log(\hat{f}_\mathbf{B}(\y_i)) + 
\mathrm{pen}(\mathbf{B}) \right\} 
\end{align}
exists and, moreover, whatever the true density $f^\star$,
\begin{align}
\label{inegOracle}
 \mathbb{E}(\dH_H^2(f^\star,\hat{f}_{\hat{\mathbf{B}}})) \leq C_{\text{oracle}} \mathbb{E}\left[ \inf_{\mathbf{B}\in \mathcal{B}_{\Lambda}} \left( \inf_{f \in F_\mathbf{B}} \KL(f^\star,f)+\mathrm{pen}(\mathbf{B})\right) \right] + 
 \frac{1\vee\tau}{n}p \log(p).
 \end{align}
\end{theorem}}

The proof is presented in \textit{Supplementary Material A}. This theorem is deduced from an adaptation of a general model selection theorem for maximum likelihood estimator developed by \cite{MassartStFlour2007}. 
This adaptation allows to focus on a random sub-collection of the whole collection of models and is also proved in \textit{Supplementary Material A}.
To use this theorem, the main assumptions to satisfy are the control of the bracketing entropy of each model in the whole collection of models and the construction of weights for each model to control 
the complexity of the collection of models. To compute the weights, we use combinatorics arguments. The control of the bracketing entropy is a classical tool to bound the Hellinger risk of the maximum 
likelihood estimator, and has already been done for Gaussian densities in \cite{Genovese2000} and \cite{Maugis2011}.

\noindent
{The assumption on the true density $f^\star$ \eqref{Bernstein2} is done because we consider a random subcollection of models $\mathcal{F}_\Lambda$ from the whole collection of models $\mathcal{F}$. 
Thanks to this assumption, we use the Bernstein inequality to control the additional randomness. The parameter $\tau$ depends on the true unknown density $f^\star$ and cannot be explicitly determined for this reason. We could do some hypothesis on the true density $f^\star$ to be able to explicit $\tau$ but we choose not to do it: e.g. under the assumption that the Kullback-Leibler divergence and the Hellinger distance are equivalent, we can explicitly determine $\tau$. Note that the parameter $\tau$ only appears in the rest term $\mathbf{r}=\frac{1\vee\tau}{n}  p \log (p)$ and not on the penalty term $\mathrm{pen}(\mathbf{B})$. Therefore, we do not need to explicit $\tau$ to select a model. Technical details are discussed in Section B.4. of the \textit{Supplementary Material A}.}

\noindent
{We remark that the Hellinger risk is upper bounded by the Kullback-Leibler divergence in \eqref{inegOracle}. For this reason, the result \eqref{inegOracle} is not exactly an oracle inequality and is called an \emph{oracle type inequality}. However, the use of the Kullback-Leibler divergence and the Hellinger distance is common for model selection theorem for Maximum Likelihood Estimator: e.g. Theorem 7.11 
in \cite{MassartStFlour2007}. Moreover, the Kullback-Leibler divergence is comparable to the Hellinger distance under some assumptions. Under these assumptions, the result \eqref{inegOracle} is exactly an oracle inequality.}

\noindent
{The collection of models \eqref{modelCollectionFmathcal} is defined such that covariance matrices have bounded eigenvalues. These bounds are useful to 
control the complexity of each model by constructing a discretization of this space. Every constant involved in \eqref{inegOracle} depends on these bounds. 
This assumption is common in non-asymptotic model selection framework. 
However, the bounds are not tractable in practice. They are calibrated in practice based on the data using the slope heuristic, as discussed in Section 4. }

To complete this analysis, we provide a second theoretical guarantee. In contrast with \cite{Lebarbier2005, Maugis2011}, we strengthen the oracle type inequality using a minimax lower bound for the risk 
between the true model and the model selected. 
Note that in Gaussian Graphical Models, lower bounds have already been obtained in other contexts \cite{Bickel2008, Cai2010}. 

\noindent
In Theorem \ref{inegalite oracle}, we have proved that we select a model as good as the oracle model in a density estimation framework. However, the bound has two extra terms: 
the penalty term $\mathrm{pen}(\mathbf{B})$ and the rest $\mathbf{r}$. The two terms give the rate of the estimator. Based on Theorem \ref{inegalite oracle} only, we do not know if the 
rate is as good 
as possible. The following theorem lower bounds the risk by a rate with the same form as the upper bound (seen as a function of $n$, $p$ and $D_\mathbf{B}$), which guarantees that we 
obtain an optimal rate.

{
\begin{theorem}
 \label{borne minimax}
 Let $\mathbf{B} \in \mathcal{B}$.
 Consider the model $F_\mathbf{B}$ defined in \eqref{myFB}, and $D_\mathbf{B}$ its dimension.
 Then, if we denote $C_{\text{minim}} = \frac{e}{4(2e+1)^2(8+\log(\lambda_M/\lambda_m))}$, for any estimator $\hat{f}_\mathbf{B}$ of $f^\star$ one has
 \begin{align}
 \label{inegMinimax}
\sup_{f^\star \in F_\mathbf{B}} \mathbb{E} (\dH_H^2(\hat{f}_\mathbf{B},f^\star)) &\geq C_{\text{minim}} \frac{D_\mathbf{B}}{n} \left(1+ \log\left(\frac{2 \lambda_M p(p-1)}{D_\mathbf{B}} \right)\right).
 \end{align}
 \end{theorem}
}

{This theorem is proved in \textit{Supplementary Material A Section C}. To obtain this lower bound, we use Birg\'e's lemma \cite{Birge2005} in conjunction with a discretization of each model, already constructed to obtain the oracle type inequality.}

\noindent
{To state Theorem \ref{borne minimax}, we assume that the parameters of the models in the collection \eqref{modelCollectionFmathcal} are bounded, which is not a strong assumption. The constant involved is explicit.}

 {Thanks to Theorems \ref{inegalite oracle} and \ref{borne minimax}, we upper bound and lower bound the Hellinger risk. The two bounds can be compared if we neglect the Kullback-Leibler term 
 (the bias term) and the rest $\mathbf{r}$. These terms are small if the collection of models is well constructed, i.e if the true density of the data is not too far from the constructed collection of models. 
 E.g., if the true model belongs to the collection of models, the two terms equal zero. 
 Thanks to Theorem \ref{borne minimax}, we say that the estimator satisfying \eqref{estimatorB} is minimax to $\mathbf{B} \in \mathcal{B}$. 
 Moreover, the lower bound is obtained for a fixed $\mathbf{B} \in \mathcal{B}$, and the rate we obtain is minimax. 
 We deduce that the performance of the procedure is as good as if we knew the structure. Consequently, our procedure is adaptive minimax, which is a strong theoretical result.}
 
{Note that the model selection procedure is optimized for density estimation and not for edge selection: the Hellinger distance and Kullback-Leibler divergence measure the differences between 
two densities from an estimation point of view. In contrast, network inference focuses on edge selection. However, we point out that the model selection procedure is only proposed in a specific context 
($n$ small), as a preliminary step (step A detailed in Section 2) prior to edge selection (step B). Although this preliminary step (step A) is not optimized for edge selection, it improves the network inference procedure as illustrated in simulated data (Section 4).}

{To conclude, let recall that these results are non-asymptotic, which means that they hold for a fixed sample size $n$, which is particularly relevant in a context where the number of observations $n$ is limited. The results are consistent with the point of view adopted in this work. }

%%%%%%%%%%%%%%%%%%%%%%%%%%%%%%%%%%%%%%%%%%%%%
%% SIMULATION
%%%%%%%%%%%%%%%%%%%%%%%%%%%%%%%%%%%%%%%%%%%%%

\section{Simulation study}
\label{Simul}

{In this section, we compare the performance of the proposed method described in Section~\ref{OurMethod} with the Cluster Graphical Lasso \cite{Tan2015} and the Graphical Lasso 
\cite{Friedman2007}. We first compare the performance of the block-diagonal covariance structure detection, i.e. step A of the proposed method (Subsection~\ref{structureDetection}), 
and then, compare the performance of the complete network inference methods (Subsection~\ref{netInf}).}

We simulate $n$ observations from a $p-$multivariate normal distribution with a null mean and a block-diagonal covariance matrix $\Sigma_\mathbf{B}$ as defined in Section~\ref{OurMethod}. 
 We fix the number of variables $p=100$, the sample size $n=70$ and the partition on variable $\mathbf{B}^\star$ with $K^\star =15$ blocks of approximately equal sizes. 
 For each block indexed by $k$, we design the $\Sigma_k$ matrix as done in \cite{GiraudHuetVerzelen2012}: $\Sigma_k= TT^t + D$ where $T$ is a random lower triangular matrix with values 
 drawn from a uniform distribution between -1 and 1, 
 and $D$ is a diagonal matrix designed to prevent $\Sigma_k$ from having eigenvalues that are too small. 
 
 {The Graphical Lasso is implemented in the R package \texttt{glasso}, version 1.7 \cite{Friedman2007}. The Cluster Graphical Lasso is based on a hierarchical clustering implemented in the R package \texttt{stats}. To detect the partition of variables (step A(2) in Section~\ref{OurMethod}), we find the connected components of the graph associated with the adjacency matrix $E_{\lambda}=[\mathbf{1}_{\{ \mid S_{j,j^\prime}\mid > \lambda\} }]_{ 1 \leq j \leq p \atop 1 \leq j^\prime \leq p }$ using a simple breadth-first search implemented in the R package \texttt{igraph} \cite{igraph}. The computation of the log-likelihood of each model (step A(3) in Section~\ref{OurMethod}) is based on the R package \texttt{mvtnorm} implementing a Cholesky decomposition \cite{mvtnorm}. To calibrate the penalty in its simplified version \eqref{kappaPP} (step A(4) in Section~\ref{OurMethod}), we use the R package \texttt{capushe} implementing the slope heuristic \cite{Baudry2012}. The complete procedure has been implemented in an R package \texttt{shock} available as \textit{Supplementary Material: the package \texttt{shock}}. }

The practical aspects of the slope heuristic are detailed in \cite{Baudry2012}: there are two methods to calibrate the penalty coefficient in \eqref{kappaPP}. One calibration method is the \textit{Slope Heuristic Dimension Jump} (SHDJ):
 the optimal coefficient $\kappa_\text{opt}$ is approximated by twice the minimal coefficient $\kappa_\text{min}$, where $\kappa_\text{min}$ corresponds to the largest dimension jump on the graph 
 representing the model dimension as a function of the coefficient $\kappa$. Another method is the \textit{Slope Heuristic Robust Regression} (SHRR): the coefficient $\kappa_\text{opt}$ is 
 approximated by twice $\kappa_\text{min}$, where $\kappa_\text{min}$ corresponds to the slope of a robust regression performed between the log-likehood and the model dimension for complex models. 
 The two methods are derived from the same heuristic and they offer two different visual checks of the adequacy of the model selection procedure to the data. They should select the same model. {Note the calibration of the more complex version of the penalty is proposed and tested in \textit{Supplementary Material B}. }The source code of the R package and the code to reproduce the simulation experiments are provided in \textit{Supplementary Materials 1, 2 and 3}.

\subsection{Block-diagonal covariance structure detection}
\label{structureDetection}

\begin{figure}[htbp!]
 \begin{center}
 \includegraphics[width= 7.7cm]{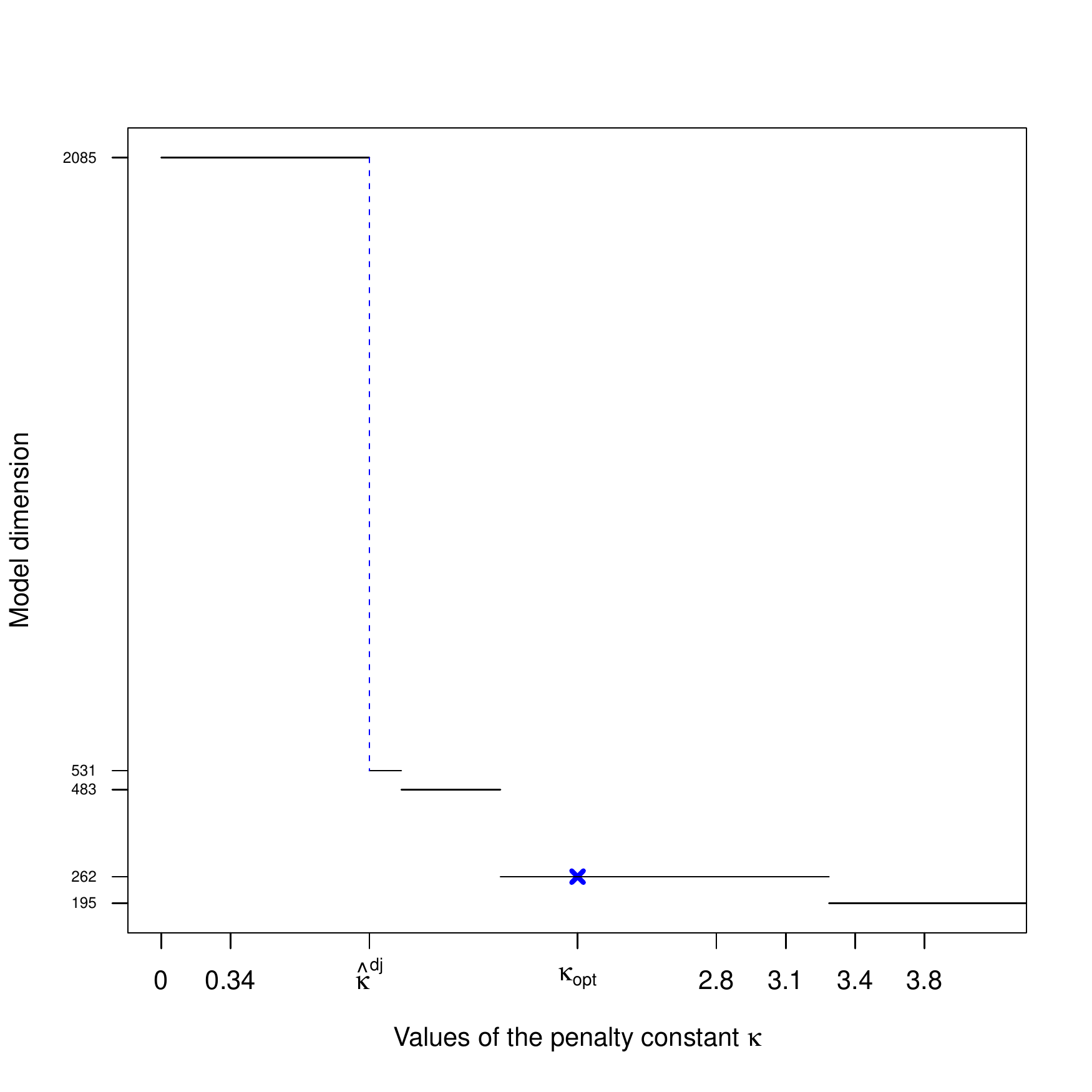}
 \includegraphics[width= 7.7cm]{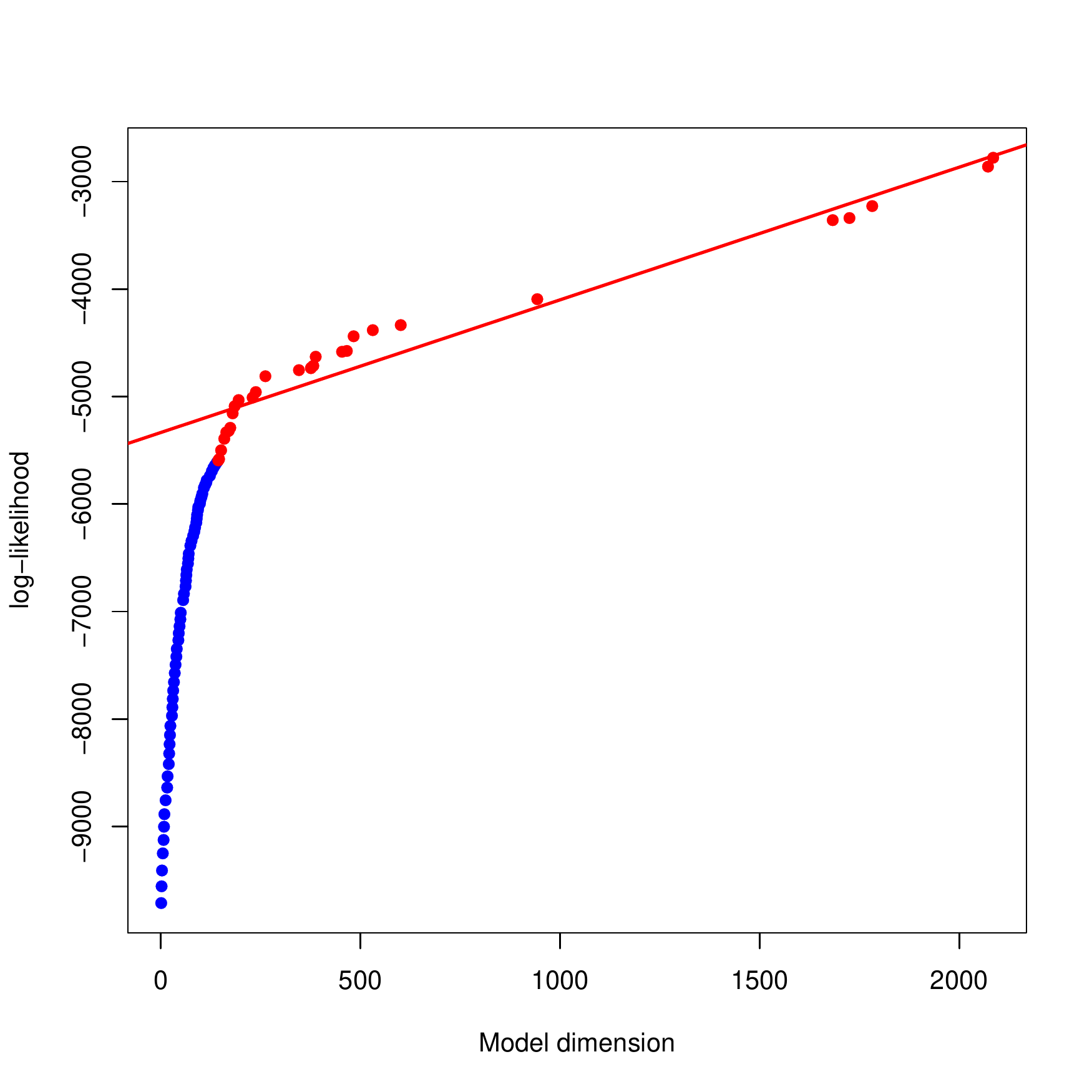}
 \end{center}
 \caption{Calibration of the $\kappa$ coefficient on a dataset simulated under a multivariate normal distribution with a block-diagonal covariance matrix $\Sigma_\mathbf{B}$ with $K^\star = 15$ blocks,
 $p=100$, $n=70$. Calibration by dimension jump (left): the dimension of the model is represented as a function of the $\kappa$ coefficient. Based on the slope heuristic, the largest jump 
 (dotted line) corresponds to the minimal coefficient $\kappa_{min}$. The optimal penalty (cross) is twice the minimal penalty. Calibration by robust regression (right): the log-likehood of 
 the model is represented as a function of the model dimension. Based on the slope heuristic, the slope of the regression (line) between the log-likehood and the model dimension for complex 
 models corresponds to the minimal coefficient $\kappa_{min}$. The optimal penalty is twice the minimal penalty. } \label{figure1}
\end{figure}

First, we investigate the ability to recover the simulated partition of variables $\mathbf{B}^\star$ based on the step A of the procedure described in Section~\ref{OurMethod}. Illustrations of the calibration of the penalty coefficient $\kappa$ are presented in Figure \ref{figure1} for one simulated dataset. The code to reproduce the simulation experiment is provided in \textit{Supplementary Material 1: Figures 1}. The largest dimension jump is easily detected on the graph representing the dimension of the model as a function of the $\kappa$ coefficient (Figure \ref{figure1} left). Likewise, we observe a linear tendency between the log-likehood and the model dimension for complex models (Figure \ref{figure1} right) and easily fit a linear regression. Both calibration methods yield the same results.

\begin{figure}[htbp!]
 \begin{center}
\includegraphics[width= 6cm]{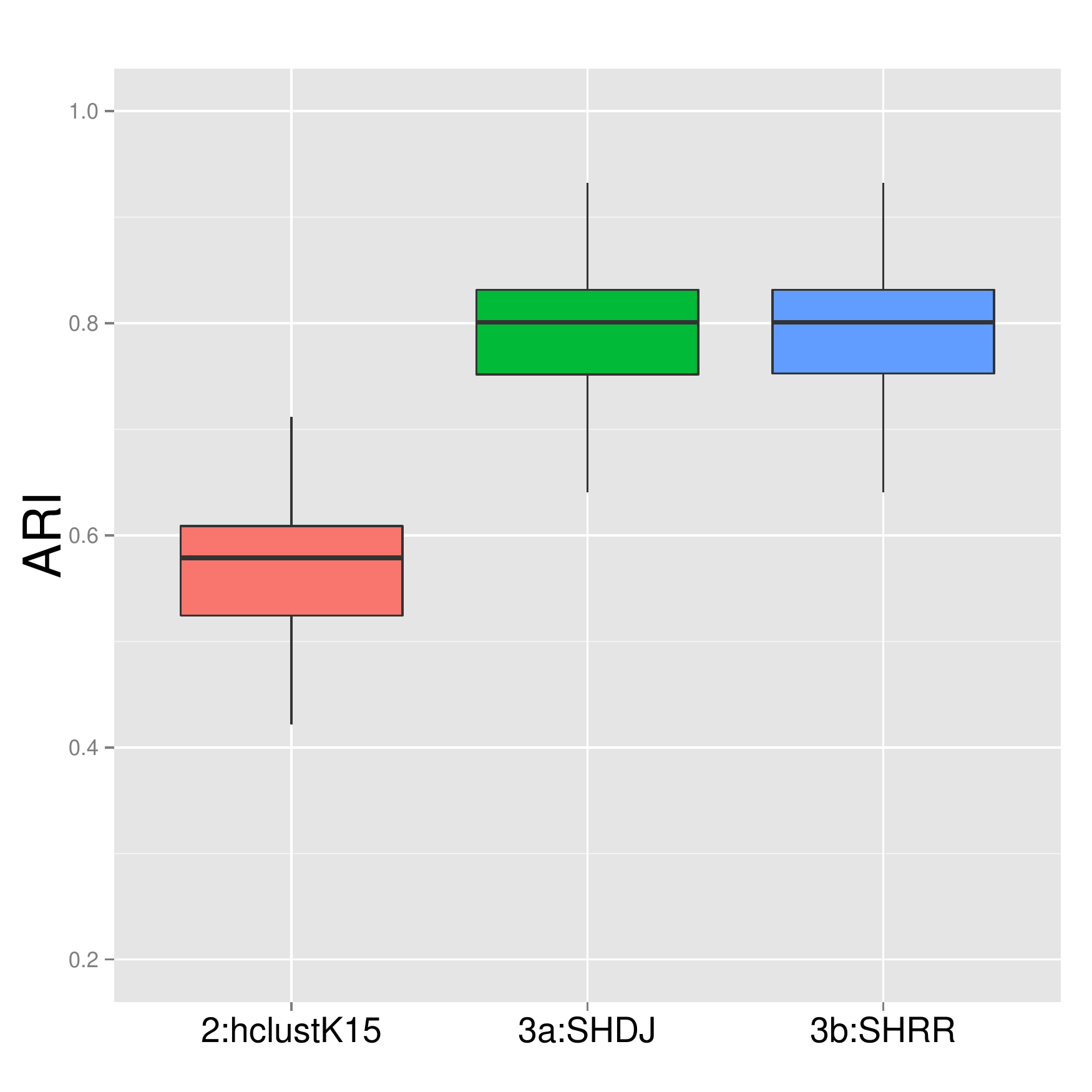}
 \end{center}
 \caption{ARI between the simulated partition and the partitions selected by slope heuristic dimension jump (SHDJ), slope heuristic robust regression (SHRR) and by average hierarchical clustering with $K=K^\star$ clusters. The ARI are computed over 100 replicated datasets simulated under a multivariate normal distribution with block-diagonal covariance matrix with $K=15$ blocks, $p=100$ variables and $n=70$ observations. } \label{figure2}
\end{figure}

In addition, we compare the partition selection methods with an average linkage hierarchical clustering with $K=K^\star$ as proposed in the cluster graphical lasso \cite{Tan2015}. Figure \ref{figure2} displays the 
{Adjusted Rand Index (ARI)} computed over 100 replicated datasets. {The ARI measures the similarity between the inferred clustering and the simulated clustering \cite{Hubert1985}. The ARI equals 1 if the two partitions match.} The code to reproduce the simulation experiment is provided in \textit{Supplementary Material 2: Figures 2 and 3}. Despite the fact that the partition with the hierarchical clustering takes as an input parameter the true number of clusters ($K=K^\star$), the ARI for the hierarchical clustering is lower than the ARI for the two slope heuristic based methods (SHRR and SHDJ) which do not need to specify the number of clusters $K$ in advance.

\subsection{Downstream network inference performance}
\label{netInf}

To illustrate the potential advantages of prior block-diagonal covariance structure detection, we compare the following strategies for network inference over 100 replicated datasets:

\begin{description}
 \item[1. Glasso:] We perform network inference using the graphical lasso on all variables, with regularization parameter $\rho$ chosen using the following $\text{BIC}^\text{net}$ criterion:
 \begin{equation}
\text{BIC}^\text{net}(\rho)=\frac{n}{2} \left(\log \det \hat{\Theta}^{(\rho)} -\text{trace}\left(S\hat{\Theta}^{(\rho)} \right)\right)- \frac{\log(n)}{2} \text{df} \hat{\Theta}^{(\rho)};
\label{BICrho2}
\end{equation}
where $\hat{\Theta}^{(\rho)}$ the solution of the graphical lasso with regularization parameter $\rho$, $S$ is the sample covariance matrix, and df the degree of freedom.
 \item[2. CGL:] We perform network inference using the cluster graphical lasso proposed in \cite{Tan2015}. First, the partition of variables is detected using an average linkage hierarchical clustering with $K=K^\star$ clusters. Note that we set the number of clusters to the true number $K^{\star}$. Subsequently, the regularization parameters in each graphical lasso problem $\rho_1, \ldots, \rho_{K^{\star}}$ are chosen from Corollary 3 of \cite{Tan2015}: the inferred network in each block must be as sparse as possible while still remaining a single connected component.
 \item[3. Inference on partitions based on model selection:] First, we detect the partition using the two variants of our non-asymptotic model selection (SHRR ou SHDJ). \hfill 
 \begin{description}
 \item[(a) SHRR:] 
 The partition $\hat{\mathbf{B}}_\text{SHRR}$ is selected using the \textit{Slope Heuristic Robust Regression}.
 \item[(b) SHDJ:] The partition $\hat{\mathbf{B}}_\text{SHDJ}$ is detected using the \textit{Slope Heuristic Dimension Jump}. 
 \hfill
 \end{description}
 Subsequently, the regularization parameters $\rho_1, \ldots, \rho_{\hat K}$ in each graphical lasso problem are chosen using the $\text{BIC}^\text{net}$ criterion:
 \begin{equation}
\text{BIC}^\text{net}(\rho_k)=\frac{n}{2} \left(\log \det \hat{\Theta}^{(\rho_k)} -\text{trace}\left({S}_{|k}\hat{\Theta}^{(\rho_k)}\right) \right)- \frac{\log(n)}{2} \text{df} \hat{\Theta}^{(\rho_k)},
\label{BICrho}
\end{equation}
where $\hat{\Theta}^{(\rho_k)}$ is the solution of the graphical lasso problem restricted to the variables in block $k$, $S_{|k}$ is the sample covariance matrix on variables belonging to the block $k$ and df the corresponding degree of freedom.
\item[4. Inference on the true partition of variables (truePart):] First, we set the partition of variables to the true partition $\mathbf{B}^\star$. Then, the regularization parameters in each graphical lasso problem $\rho_1, \ldots, \rho_{K^{\star}}$ are chosen using the $\text{BIC}^\text{net}$ criterion \eqref{BICrho}.
 \end{description} 

We compare the performance of the five methods using the sensitivity ($\mathit{Sensitivity} = \mathit{TP} / (\mathit{TP} + \mathit{FN})$), the specificity ($\mathit{Specificity} = \mathit{TN} / (\mathit{TN}+\mathit{FP})$) and the False Discovery Rate (FDR) ($\mathit{FDR} = \mathit{FP} / (TP + FP)$) where $TN, TP, FN, FP$ are respectively the number of true negative, true positive, false negative, false positive dependencies detected. A network inference procedure is a compromise between sensitivity and specificity: we are looking for a high sensitivity, which measures the proportion of dependencies (presence of edges) that are correctly identified, and a high specificity, which measures the proportion of independencies (absence of edges) that are correctly identified. The False Discovery Rate is the proportion of dependencies wrongly detected. The value of sensitivity, specificity and False Discovery Rate are computed over 100 replicated datasets as illustrated in Figure \ref{figure3}. The code to reproduce the simulation experiment is provided in \textit{Supplementary Material 2: Figures 2 and 3}. As expected, the true partition strategy (truePart) performs the best: based on the true partition of variables, the network inference problem is easier because we solve problems of smaller dimension. The proposed strategies, based on the SHRR and SHDJ partitions, improve network inference compared to a simple graphical lasso on the set of all variables (glasso) or compared to the cluster graphical lasso (CGL).

 \begin{figure}[htbp!]
 \begin{center}
\includegraphics[width=5.2cm]{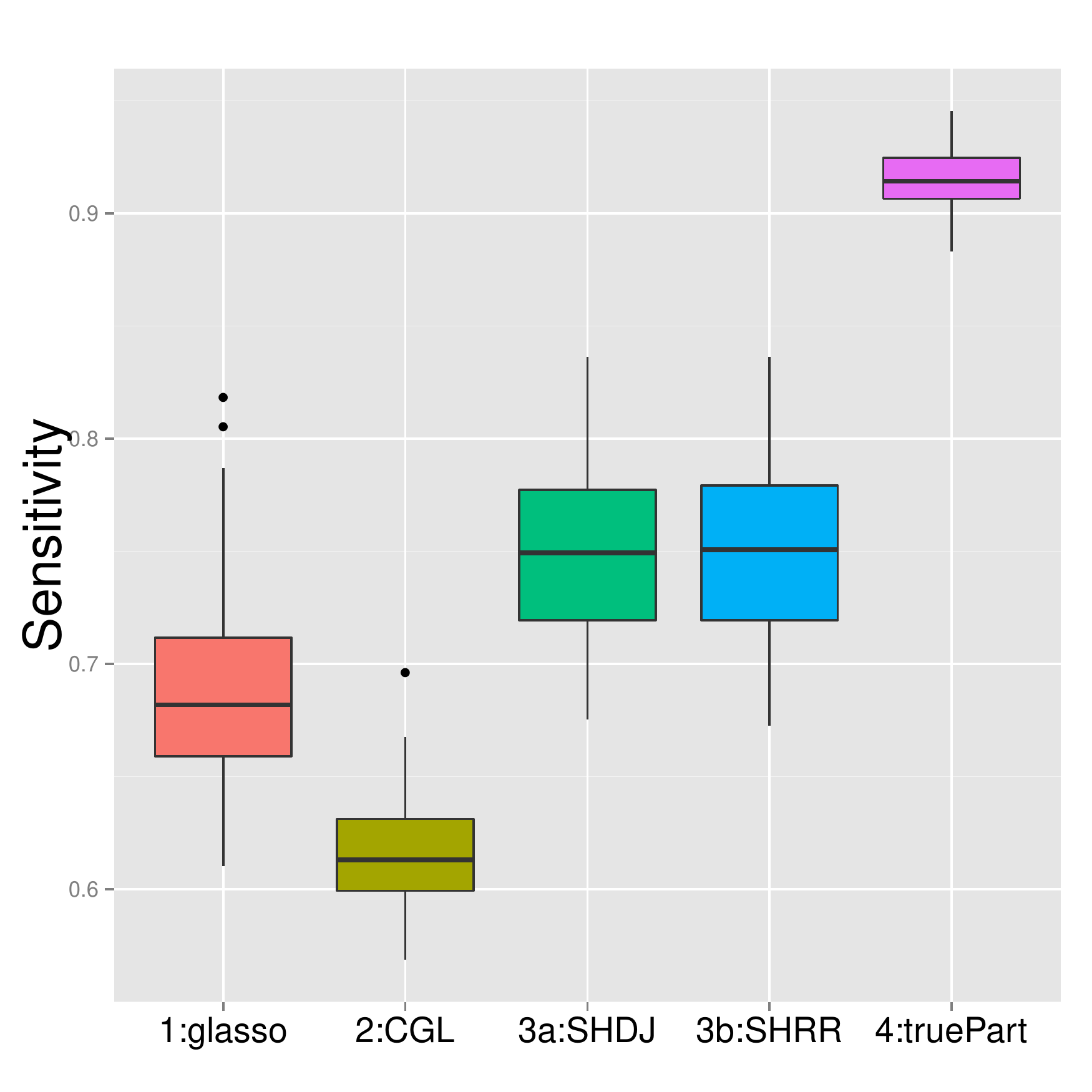}
\includegraphics[width=5.2cm]{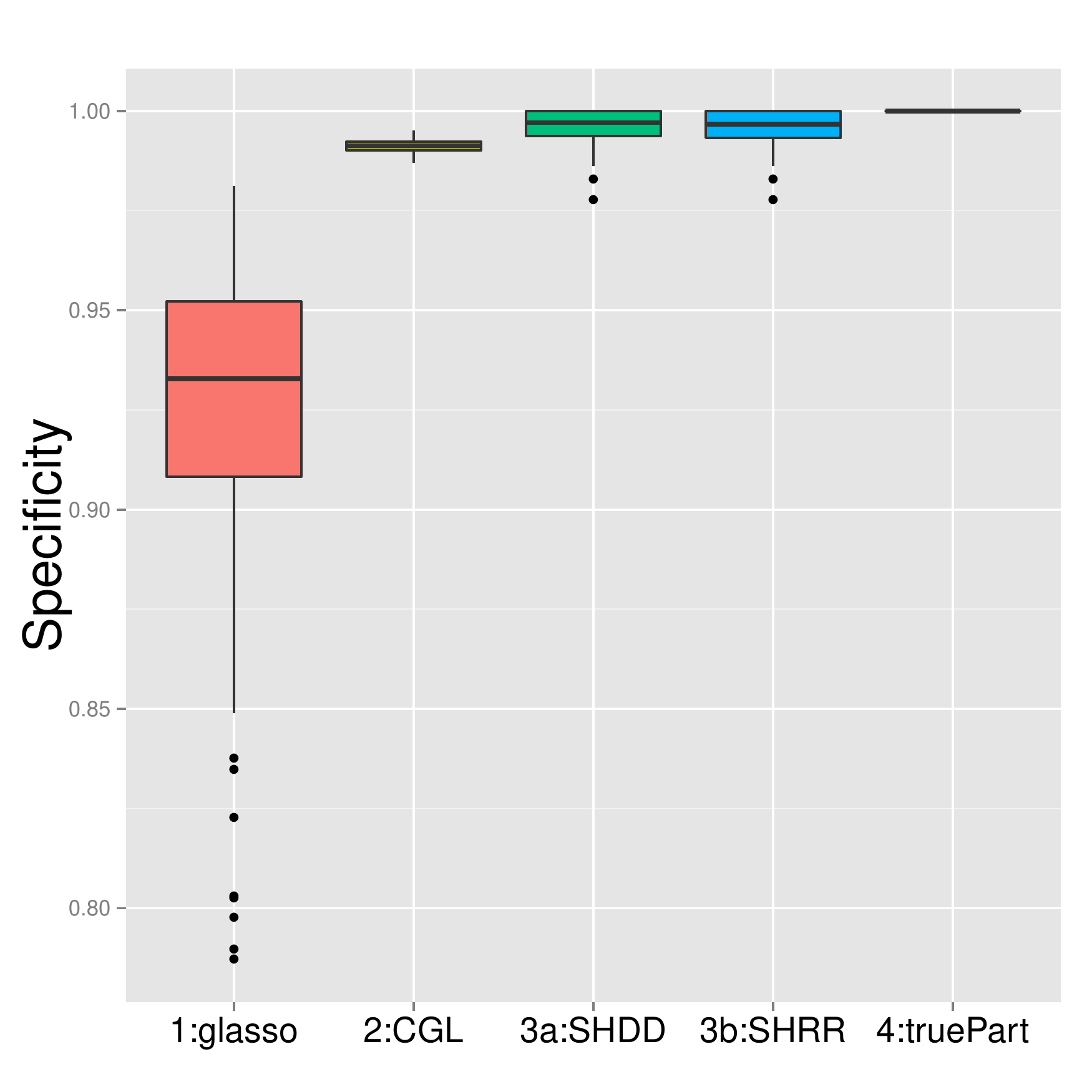}
\includegraphics[width=5.2cm]{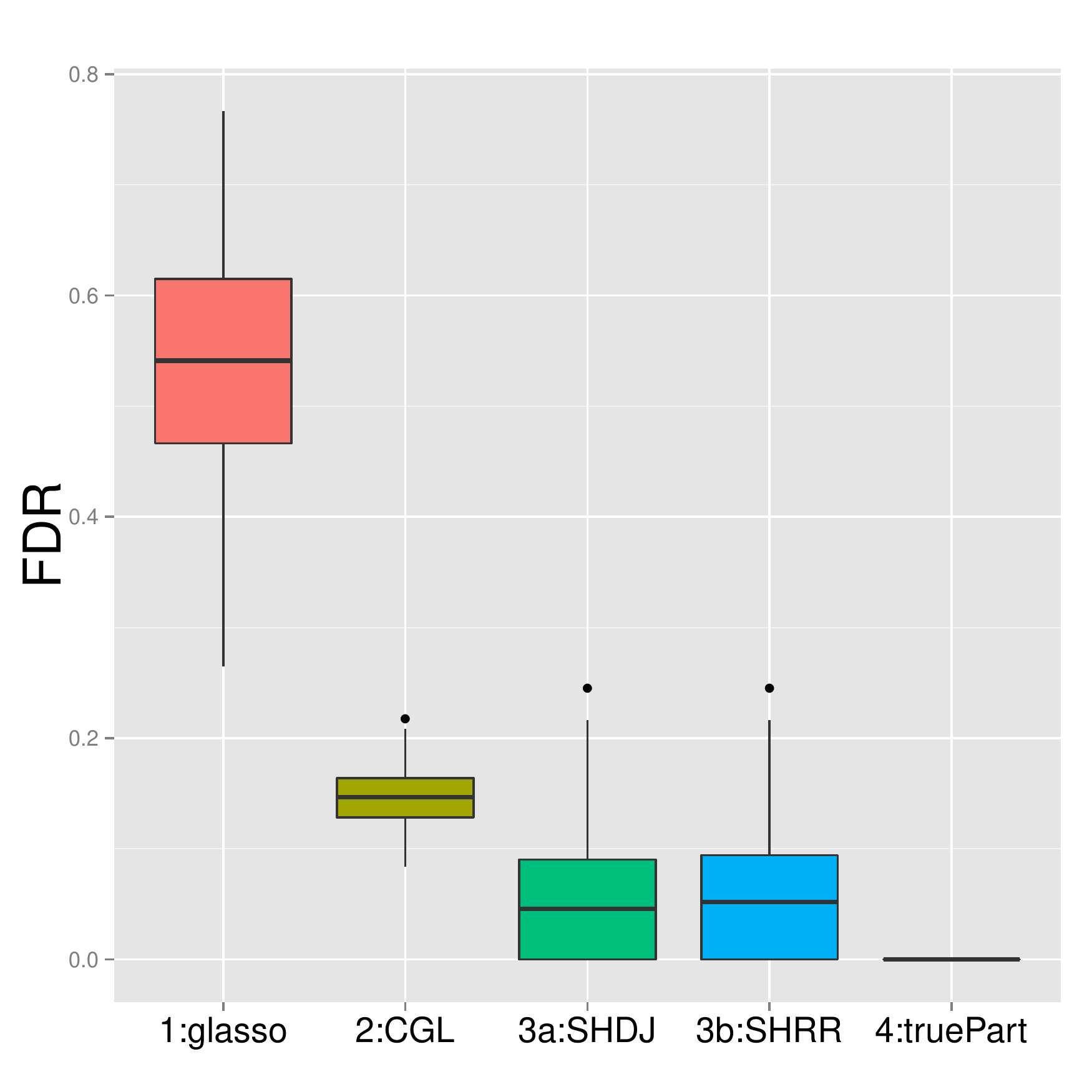}
 \end{center}
 \caption{ Performance of network inference methods (glasso: graphical lasso on the set of all variables, CGL: cluster graphical lasso, BIC: network inference based on the partition of variables $\hat{\mathbf{B}}_{\text{BIC}}$, SSHR: network inference based on the partition of variables $\hat{\mathbf{B}}_{\text{SHRR}}$, SHDJ: network inference based on the partition of variables $\hat{\mathbf{B}}_{\text{SHDJ}}$ and truePart: network inference based on the partition of variables $\mathbf{B}^\star$) measured by the sensitivity, the specificity and the False Discovery Rate (FDR) of the inferred graph over 100 replicated datasets simulated under a $p-$multivariate normal distribution with a null mean $\mathbf{0}$ and a block-diagonal covariance matrix $\Sigma_{\mathbf{B}^\star}$ with $p=100$, $K=15$, $n=70$ and clusters of approximately equal sizes. } 
 \label{figure3}
\end{figure}

%%%%%%%%%%%%%%%%%%%%%%%%%%%%%%%%%%%%%%%%%%%%%
%% REAL DATA ANALYSIS
%%%%%%%%%%%%%%%%%%%%%%%%%%%%%%%%%%%%%%%%%%%%%
\newpage
\section{Real data analysis}
\label{RealData}

Pickrell {\it et al.} analyzed transcriptome expression variation from 69 lymphoblastoid cell lines derived from unrelated Nigerian individuals \cite{Pickrell2010}. The expression of 52580 genes across 69 observations was measured using RNA-seq. The data is extracted from the Recount database \cite{Frazee2011}. After filtering weakly expressed genes using the \texttt{HTSFilter} package \cite{Rau2013}, we identified the $200$ most variable genes among the 9191 remaining genes, and restrict our attention to this set of genes for the following network inference analysis. The code to reproduce the analysis is provided in \textit{Supplementary Material 3}. 

\begin{figure}[htbp!]
 \begin{center}
 \includegraphics[width= 8cm]{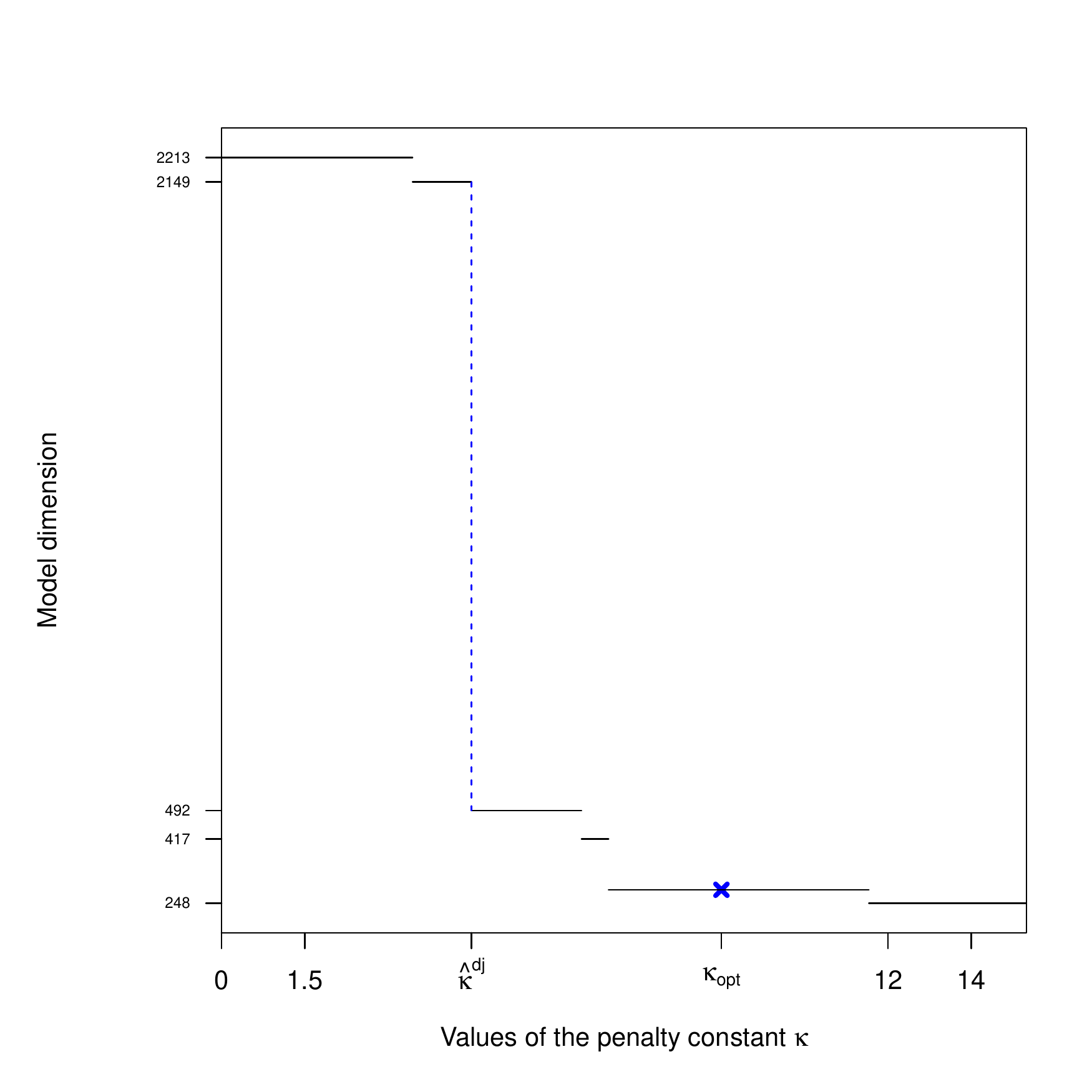}
 \includegraphics[width= 8cm]{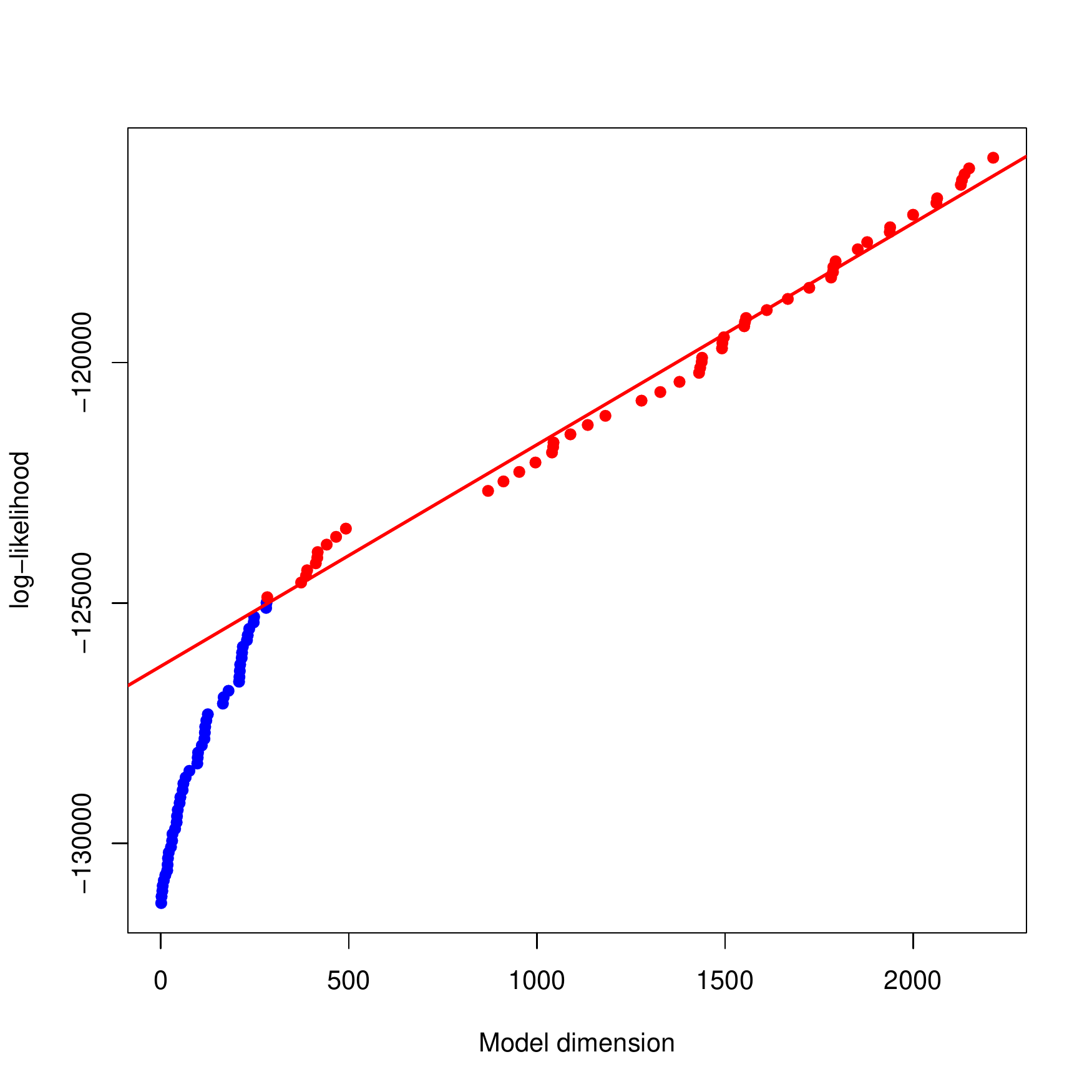}
 \end{center}
 \caption{Calibration of the $\kappa$ coefficient on the $200$ most variable genes extracted from the \cite{Pickrell2010} dataset. Calibration by robust regression (left) and by dimension jump (right).
 In both cases, the optimal penalty is twice the minimal penalty.}
\label{logPickrell}
\end{figure}

 First, we select the partition $\hat{\mathbf{B}}$ using model selection as described in equation \eqref{kappaPP}. The log-likelihood increases with the number of parameters to be estimated in the model as displayed in Figure \ref{logPickrell}. We notice a linear tendency in the relationship between the log-likelihood and the model dimension for complex models (points corresponding to a model dimension higher than 500).
This suggests that the use of the slope heuristic is appropriate for selecting a partition $\hat{\mathbf{B}}$. The model selected by SHDJ and by SHRR described in Section~\ref{OurMethod} are the same. The number of blocks detected is $\hat K_\text{SH}=150$ and the corresponding model dimension is $D_{\hat{\mathbf{B}}_\text{SH}}=283$. The partition $\hat{\mathbf{B}}_\text{SH}$ yields 4 blocks of size $18,13,8$ and $5$, 4 blocks of size 3, 2 blocks of size 2 and 140 blocks of size 1. The partition selected by the slope heuristic offers a drastic reduction of the number of parameters to infer, as compared with the graphical lasso performed on the full set of variables, which corresponds to a total of $D = 19900$ parameters to estimate.

\begin{figure}[htbp!]
 \begin{center}
 \includegraphics[width= 3.8cm]{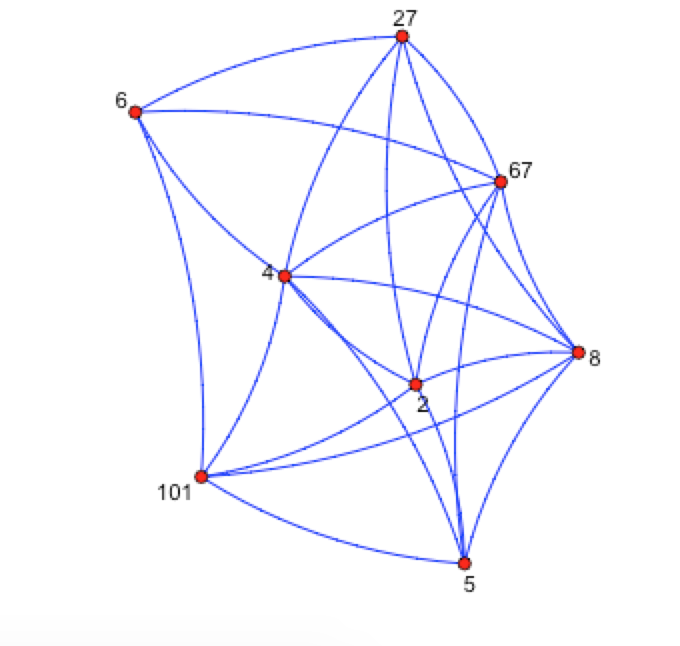} 
 \includegraphics[width= 3.8cm]{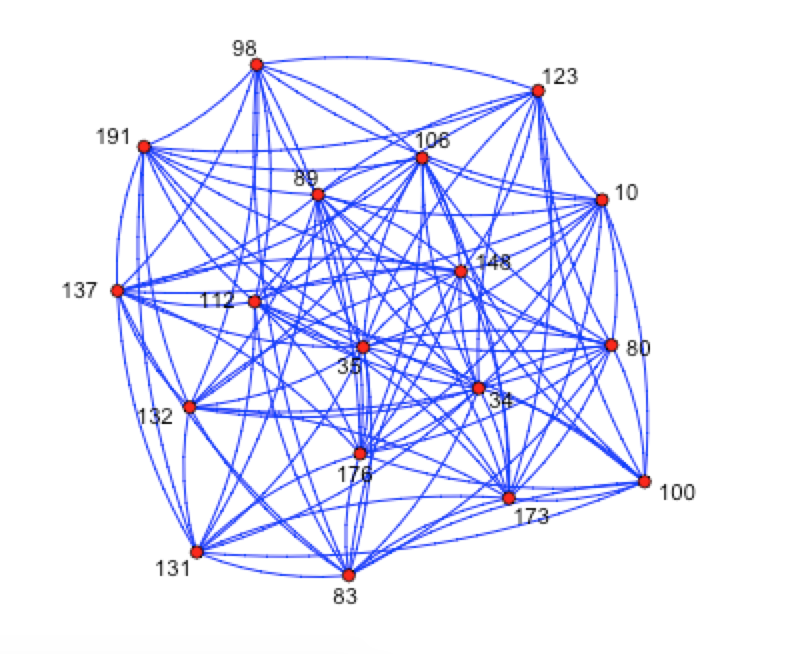}
 \includegraphics[width= 3.8cm]{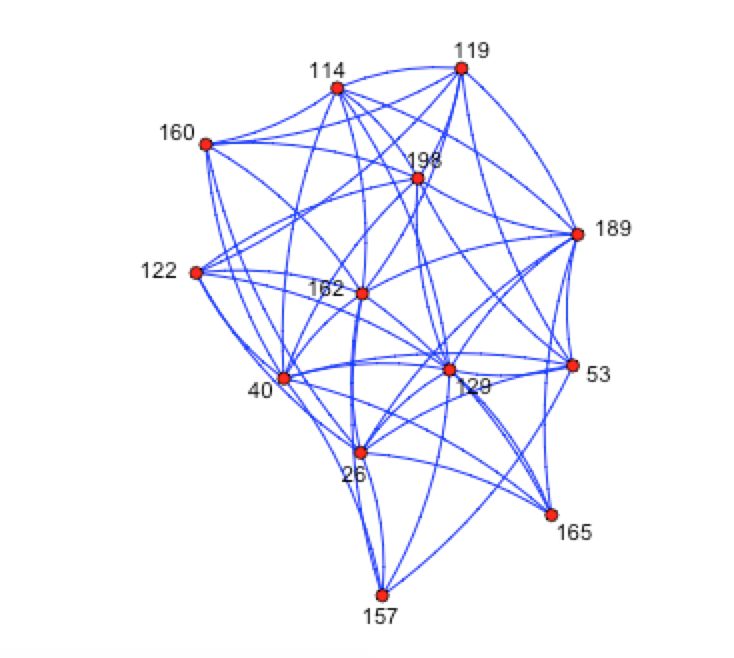}
 \includegraphics[width= 3.8cm]{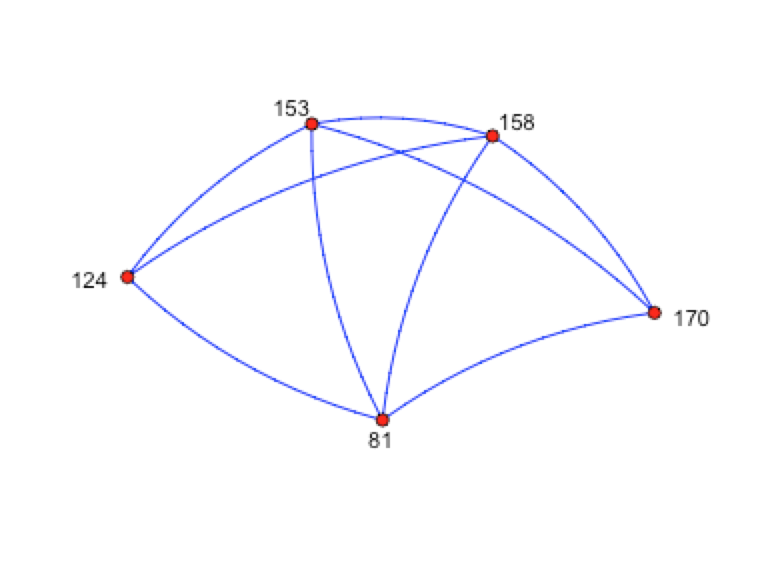}
 \end{center}
 \caption{Networks inferred on the four largest components detected by slope heuristic. Regularization parameters in each set of variables are chosen using the $\text{BIC}^\text{net}$ criterion \eqref{BICrho}. Numbers indicate gene labels. } \label{networkDDSE}
\end{figure}

%\paragraph{Network inference.} 
The networks within each cluster of variables are inferred using the graphical lasso algorithm of Friedman \cite{Friedman2007} implemented in the \texttt{glasso} package, version 1.7.
 The regularization parameter for the graphical lasso on the set of all variables is chosen using the $\text{BIC}^\text{net}$ criterion \eqref{BICrho2}. The model inferred based on partition $\hat{\mathbf{B}}_\text{SH}$ is more parsimonious and easier to interpret than the model inferred on the full set of variables. An illustration of inferred networks in the four largest connected components of the partition $\hat{\mathbf{B}}_\text{SH}$ are displayed on Figure \ref{networkDDSE}. These four networks might be good candidates for further study.

%%%%%%%%%%%%%%%%%%%%%%%%%%%%%%%%%%%%%%%%%%%%%
%% DISCUSSION
%%%%%%%%%%%%%%%%%%%%%%%%%%%%%%%%%%%%%%%%%%%%%

\section{Discussion}

In this paper, we propose a non-asymptotic procedure to detect a block diagonal structure for covariance matrices in GGMs. Our non-asymptotic approach is supported by theoretical results:
 an oracle type inequality ensures that the model selected based on a penalized criterion is close to the oracle, {\it i.e.} the best model among our family. Moreover, we obtain a minimax lower bound 
 of the risk between the true model and the model selected among the collection of models, which ensures that the model selection procedure is optimal, {\it i.e.} adaptive minimax to the block-structure.

The method we propose is easy to implement in practice and fast to compute. The calibration of the $\kappa$ coefficient by robust regression and dimension jump encounters no particular difficulty. Moreover, graphical representations of the two slope heuristic calibration methods allow to visualize the calibration performance.

Note that our non-asymptotic results hold for a fixed number of sample, which is a typical case in many real applications. Although GGMs are widely used in practice, limited sample sizes typically
force the user to restrict the number of variables. Usually, this restriction is performed manually based on prior knowledge on the role of variables. Here, our procedure allows to select relevant
subsets of variables based on a data-driven criterion. This procedure is of great practical interest to estimate parameters in GGMs when the sample size is small. We apply it on a high-throughput genomic dataset
but the procedure may be useful for other types of data with low sample size (e.g., neuroscience, sociology).

\section{Appendix: Theoretical proofs}
In this \textit{Supplementary Material A}, we give details about the proofs of the theoretical results presented in the paper.

First, we recall the collection of models as defined in the article.
In Section \ref{Discretization}, we describe a discretization of the collection of models, which is useful to prove the oracle type inequality and the lower bound.
In Section \ref{oracle type inequality proof}, we prove the oracle type inequality. %To do so, we first state and prove a model selection theorem for maximum likelihood estimator when we focus on a random subcollection of the whole collection of models.
% 
% we prove the oracle type inequality (Theorems 3.1) and the minimax lower bound (Theorem 3.2) presented in the manuscrit. 
% First, we describe a discretization of the collection of models used, which is useful for both proofs.
% In Section , we prove the oracle type inequality. To do so, we first state and prove Theorem \ref{ThmgeneralModelSelectionTheorem}, a generalization of Theorem 7.11 from \cite{MassartStFlour2007}, a model selection theorem for the maximum likelihood estimator (MLE) to a random sub-collection of models.
%Subsequently, we prove that our collection of models satisfies all the assumptions of this theorem, and then deduce the oracle type inequality.
In Section \ref{lower bound proof}, we prove the minimax lower bound using Birg\'e's lemma.%) in conjunction with the discretization of the collection of models obtained in Section \ref{Discretization}. 

We recall the main notations as defined in the article.
Let $(\y_1,\ldots,\y_n)$ be a sample of $\mathbb{R}^p$ of size $n$. We consider the following collection of models, for $\mathbf{B} \in \mathcal{B}$ the partition of variables into $K$ blocks,
\begin{align}
 F_{\mathbf{B}} &= \left\{ f_{\mathbf{B}} =\phi _p(0,{{\SSigma}}_{\mathbf{B}}) \text{ with } \SSigma_{\mathbf{B}} \in S_{\mathbf{B}} \right\}, \\
S_{\mathbf{B}} &= \left\{ \SSigma_{\mathbf{B}} \in \mathbb{S}_p^{++}(\mathbb{R})\left| 
 \begin{array}{l} \lambda_m \leq {\min}(\text{sp}({\SSigma}_{\mathbf{B}})) \leq {\max}(\text{sp}({\SSigma}_{\mathbf{B}})) \leq \lambda_M, \\
%  \lambda_m \leq \Lambda_{\min} (\SSigma_{\mathbf{B}}) \leq \Lambda_{\max}(\SSigma_{\mathbf{B}}) \leq \lambda_M,\\
    \SSigma_{\mathbf{B}} = P_\sigma 
%  \resizebox{3cm}{!}{
 \begin{pmatrix}
  {\SSigma}_1 & 0 & 0\\
  0 & \ddots & 0 \\
  0 & 0 & {\SSigma}_K
 \end{pmatrix}%}
 P_\sigma^{-1}, \\
  {\SSigma}_k \in \mathbb{S}_{p_k}^{++}(\mathbb{R}) \text{ for } k \in \{1,\ldots,K\} \\
 \end{array}
 \right.\right\},
\label{myFB2}
\end{align}
where $\mathbb{S}_p^{++}(\mathbb{R})$ is the set of positive semidefinite matrices, $\lambda_m$ and $\lambda_M$ are real numbers, 
${\min}(\text{sp}({\SSigma}_{\mathbf{B}})), {\max}(\text{sp}({\SSigma}_{\mathbf{B}}))$ are the smallest and highest eigenvalues of $\SSigma_{\mathbf{B}}$ and $P_\sigma$ is a permutation matrix. 
We also note $e_m$ and $e_M$ 
the smallest and the largest values of $\SSigma_{\mathbf{B}}$ (remark that $e_m$ and $e_M$ could be bounded by $\lambda_M$, but we use the notations $e_m$ and $e_M$ to simplify the reading).

\subsection{Model collection and discretization}
\label{Discretization}
The aim of this section is to discretized the set $S_{\mathbf{B}}$.
Our theoretical results rely on combinatorial arguments, such as an assumption on the bracketing entropy for the oracle type inequality and the Birg\'e's lemma for the lower bound. 
Then, we develop in this section 
in details the discretization used later on.

We distinguish between low- and high-dimensional models. 
For low-dimensional models, we construct  a thiner discretization, whereas in the high-dimensional models, the densities are closer and the situation is more complicate.
\subsubsection{Discretization of the set of adjacency matrices for low-dimensional models}

We denote by $Adj(A)$ the adjacency matrix associated to a covariance matrix $A$.
Let $\mathbf{B}=(\mathbf{B}_1,\ldots,\mathbf{B}_K) \in \mathcal{B}$.
For a given matrix $\SSigma_{\mathbf{B}} \in S_{\mathbf{B}}$, we may identify a corresponding adjacency matrix $A_{\mathbf{B}} = Adj(\SSigma_{\mathbf{B}})$. This matrix of size $p^2$ could be summarized by the vector of concatenated upper triangular vectors.
Then, we construct a discrete space for $\{0,1\}^{p(p-1)/2}$ which is in bijection with
\begin{align*}
 \mathcal{A}_{\mathbf{B}} &= \left\{ A_{\mathbf{B}} \in \mathbb{S}_p(\{0,1\}) | \exists \SSigma_{\mathbf{B}} \in S_{\mathbf{B}} \text{ s.t. } Adj(\SSigma_{\mathbf{B}}) = A_{\mathbf{B}} \right\}.
\end{align*}
First, we focus on the set $\{0,1\}^{p(p-1)/2}$.
\begin{lemma}
\label{constructionTheta}
 Let $\{0,1\}^{p(p-1)/2}$ be equipped with Hamming distance $\delta$.
 Let $\{0,1\}^{p(p-1)/2}_{\mathbf{B}}$ be the subset of $\{0,1\}^{p(p-1)/2}$ of vectors for which the corresponding graph has structure $\mathbf{B}$.

For every $\alpha \in (0,1)$, let $\beta \in(0,1)$ such that $D_{\mathbf{B}} \leq \alpha\beta p(p-1)/2$. There exists some subset $\mathcal{R}(\alpha)$ of $\{0,1\}^{p(p-1)/2}_{\mathbf{B}}$ 
with the following properties:
\begin{align}
\delta(r,r^{'}) &> 2(1-\alpha) D_{\mathbf{B}} \text{ for every } (r,r^{'}) \in \mathcal{R}(\alpha)^2 \text{ with } r \neq r^{'}; \label{1}\\
\log |\mathcal{R}(\alpha)| &\geq \rho D_{\mathbf{B}} \log \frac{p(p-1)}{2 D_{\mathbf{B}}} -p\left(\frac{1}{2} + \log(p)\right);
\end{align}
where $\rho=-\alpha(-\log(\beta)+\beta-1)/\log(\alpha \beta)$ and $D_{\mathbf{B}}=\sum_{1\leq k \leq K} p_k(p_k-1)/2$.
\end{lemma}
\begin{demo}
 Let $\mathcal{R}$ be a maximal subset of $\{0,1\}^{p(p-1)/2}_{\mathbf{B}}$ satisfying property \eqref{1}.
 Then the closed balls with radius $\epsilon$ whose belongs to $\mathcal{R}$ cover $\{0,1\}^{p(p-1)/2}_{\mathbf{B}}$.
 We remark that $x \mapsto P_\sigma x P_\sigma^{-1}$ is a group action, isometric and transitive on $\{0,1\}^{p(p-1)/2}_{\mathbf{B}}$.
 
 Hence,
 $$| \{0,1\}^{p(p-1)/2}_{\mathbf{B}}| \leq \sum_{x \in \mathcal{R}} |\mathscr{B}_{\{0,1\}^{p(p-1)/2}_{\mathbf{B}}}(x,\epsilon)| = |\mathcal{R}| |\mathscr{B}_{\{0,1\}^{p(p-1)/2}_{\mathbf{B}}}(x^0,\epsilon)| $$
 for every $x^0 \in \mathcal{R}$, where $\mathscr{B}_A(x,r) = \{ y \in A | \delta(x,y) \leq r\}$.

 Our proof is similar to the proof of Lemma 4.10 in \cite{MassartStFlour2007}. We consider:
 $$[\{0,1\}^{p(p-1)/2}]_D= \left\{ x \in\{0,1\}^{p(p-1)/2} | \delta(0,x)=D\right\}.$$
 
 Let $x^0 \in \mathcal{R}$.
 Let $\alpha \in (0,1), \beta \in (0,1)$ such that $D \leq \alpha \beta p(p-1)/2$.
 Due to \cite{MassartStFlour2007}, we prove that
 \begin{align*}
 |\mathscr{B}_{[\{0,1\}^{p(p-1)/2}]_D}(x^0, 2(1-\alpha)D)| &\leq \frac{\exp\left(- \rho D\log(p(p-1)/2{D})\right)}{\dbinom{p(p-1)/2}{{D}}}.
 \end{align*}
 with $\rho = - \alpha (- \log(\beta)+\beta-1)/\log(\alpha \beta)$. It relies on exponential hypergeometric tail bounds.
 %where $h(u)=(1+u)\log(1+u)-u$.
 
 Nevertheless, as $\{0,1\}^{p(p-1)/2}_{\mathbf{B}} \subset [\{0,1\}^{p(p-1)/2}]_{D_{\mathbf{B}}}$, for $D_{\mathbf{B}} = \sum_{k=1}^K p_k(p_k-1)/2$,
 $$|\{0,1\}^{p(p-1)/2}_{\mathbf{B}} | \leq |\mathcal{R}|\frac{\exp\left(- \rho D_{\mathbf{B}}\log(p(p-1)/2{D_{\mathbf{B}}})\right)}{\dbinom{p(p-1)/2}{{D_{\mathbf{B}}}}}.$$
 
As 
 $$|\{0,1\}^{p(p-1)/2}_{\mathbf{B}} | \geq 1$$
and
\begin{align*}
 \dbinom{p(p-1)/2}{{D_{\mathbf{B}}}}\geq 1,
 \end{align*}
 we get
 \begin{align*}
 |\mathcal{R}(\alpha)| &\geq \exp \left(\rho D_{\mathbf{B}}\log(p(p-1)/2{D_{\mathbf{B}}})\right).
 \end{align*}
\end{demo}

\subsubsection{Discretization of the set of covariance matrices for low-dimensional models}
In this section, we use Lemma \ref{constructionTheta} to deduce a discretization of the set of covariance matrices for low-dimensional models.
\begin{proposition}
\label{constructionSBDisc}
Let $\alpha \in (0,1)$ and $\beta \in (0,1)$ such that $D_{\mathbf{B}} \leq \alpha \beta p(p-1)/2$.
 Let $\mathcal{R}(\alpha)$ as constructed in Lemma \ref{constructionTheta}, and its equivalent $\mathcal{A}_{\mathbf{B}}^{\text{disc}}(\alpha)$ for adjacency matrices.
 Let $\epsilon >0$.
 Let:
 $$S_{\mathbf{B}}^\text{disc}(\epsilon,\alpha) = \left\{ \SSigma \in \mathbb{S}_p^{++}(\mathbb{R}) | \text{Adj}(\SSigma) \in \mathcal{A}_{\mathbf{B}}^{\text{disc}}(\alpha),
 \SSigma_{i,j} = \sigma_{i,j} \epsilon, \sigma_{i,j} \in \left[\frac{-\lambda_M}{\epsilon},\frac{\lambda_M}{\epsilon}\right] \cap \mathbb{Z} \right\}.$$
 Then,
 \begin{align*}
 ||\SSigma-\SSigma'||_2^2 &\geq2(1-\alpha)D_{\mathbf{B}} \wedge \epsilon \text{ for every } (\SSigma, \SSigma') \in (S_{\mathbf{B}}^\text{disc}(\epsilon,\alpha))^2 \text{ with } \SSigma \neq \SSigma';\\
 \log | S_{\mathbf{B}}^\text{disc}(\epsilon,\alpha) |&\geq \rho D_{\mathbf{B}} \log \left( \left\lfloor\frac{2 \lambda_M}{\epsilon}\right\rfloor \frac{p(p-1)}{2 D_{\mathbf{B}}} \right).
 \end{align*}

\end{proposition}
\begin{demo}
 Let $(\SSigma, \SSigma') \in (S_{\mathbf{B}}^\text{disc}(\epsilon,\alpha))^2$ with $\SSigma \neq \SSigma'$.
 If $\SSigma$ and $\SSigma'$ are close, either they have the same adjacency matrix and they differ only on a coefficient or they differ in their adjacency matrices.
 In the first case, $||\SSigma-\SSigma'||_2^2 \geq \epsilon$.
 In the second case, $||\SSigma-\SSigma'||_2^2 \geq 2(1-\alpha) D_{\mathbf{B}}$.
 Then,
 $$||\SSigma-\SSigma'||_2^2 \geq2(1-\alpha)D_{\mathbf{B}} \wedge \epsilon,$$
 this minimum depending on $\alpha$ and $\epsilon$.
 
\end{demo}

\begin{corollary}
\label{constructionSBDisc2}
 If $D_{\mathbf{B}} \leq p(p-1)/8$, let $\mathcal{R}(3/4)$ as constructed in Lemma \ref{constructionTheta}, and its equivalent $\mathcal{A}_{\mathbf{B}}^{\text{disc}}(3/4)$ for adjacency matrices.
 Let $\epsilon >0$, and
 $$\tilde{S}_{\mathbf{B}}^\text{disc}(\epsilon) = \left\{ \SSigma \in \mathbb{S}_p^{++}(\mathbb{R}) | \text{Adj}(\SSigma) \in \mathcal{A}_{\mathbf{B}}^{\text{disc}}(3/4),
 \SSigma_{i,j} = \sigma_{i,j} \epsilon, \sigma_{i,j} \in \left[\frac{-\lambda_M}{\epsilon},\frac{\lambda_M}{\epsilon}\right] \cap \mathbb{Z} \right\}.$$
 Then,
 \begin{align*}
 ||\SSigma-\SSigma'||_2^2 &\geq \frac{1}{2} D_{\mathbf{B}} \wedge \epsilon \text{ for every } (\SSigma, \SSigma') \in (\tilde{S}_{\mathbf{B}}^\text{disc}(\epsilon))^2 \text{ with } \SSigma \neq \SSigma'\\
 \log | \tilde{S}_{\mathbf{B}}^\text{disc}(\epsilon) |&\geq \rho D_{\mathbf{B}} \log \left( \left\lfloor\frac{2 \lambda_M}{\epsilon}\right\rfloor \frac{p(p-1)}{2 D_{\mathbf{B}}} \right)-p\left(\frac{1}{2} + \log(p)\right).
 \end{align*}
 with $\rho \geq 0.233$.

\end{corollary}

\subsubsection{Discretization of the set of covariance matrices for high-dimensional models}
In this section, we deal with high-dimensional models. The following lemma and its proof could be found in \cite{MassartStFlour2007}, Lemma 4.7.
\begin{lemma}[Varshamov-Gilbert]
\label{Varshamov}
 Let $\{0,1\}^{D_{\mathbf{B}}}$ be equipped with Hamming distance $\delta$. Given $ \alpha \in (0,1)$, there exists some subset $\Theta(\alpha)$ of $\{0,1\}^{D_{\mathbf{B}}}$ 
 with the following properties:
  \begin{align*}
   \delta(r,r') &> (1-\alpha)\frac{D_{\mathbf{B}}}{2} \text{ for every } (r,r')\in \Theta(\alpha)^2 \text{ with } r\neq r' \\
   \log |\Theta(\alpha)| &\geq \frac{\rho D_{\mathbf{B}}}{2}
  \end{align*}
where $\rho = (1+\alpha) \log(1+\alpha) + (1-\alpha) \log(1-\alpha)$. 
\end{lemma}
We deduce a discretization of the set of covariance matrices for high-dimensional models.
\begin{proposition}
\label{highDimModel}
 If $D_{\mathbf{B}} \geq p(p-1)/8$, let $\Theta(1/2)$ as constructed in Lemma \ref{Varshamov}, and its equivalent $\check{\mathcal{A}}^{\text{disc}}(1/2)$ for adjacency matrices. Let $\epsilon >0$, and let
 $$\check{S}_{\mathbf{B}}^\text{disc}(\epsilon) = \left\{ \SSigma \in \mathbb{S}_p^{++}(\mathbb{R}) | \text{Adj}(\SSigma) \in \check{\mathcal{A}}^{\text{disc}}(1/2),
 \SSigma_{i,j} = \sigma_{i,j} \epsilon, \sigma_{i,j} \in \left[\frac{-\lambda_M}{\epsilon},\frac{\lambda_M}{\epsilon}\right] \cap \mathbb{Z} \right\}.$$
 Then,
 \begin{align*}
 ||\SSigma-\SSigma'||_2^2 &\geq \frac{1}{2} D_{\mathbf{B}} \wedge \epsilon \text{ for every } (\SSigma, \SSigma') \in (\check{S}_{\mathbf{B}}^\text{disc}(\epsilon))^2 \text{ with } \SSigma \neq \SSigma'\\
 \log | \check{S}_{\mathbf{B}}^\text{disc}(\epsilon) |&\geq \frac{ D_{\mathbf{B}}}{4} \log \left( \left\lfloor\frac{2 \lambda_M}{\epsilon}\right\rfloor \right).
 \end{align*}
\end{proposition}
\begin{demo}
 As $\{0,1\}^{D_{\mathbf{B}}} \subset \{0,1\}^{p(p-1)/2}_D$, we use the Varshamov-Gilbert's lemma.
 With $\alpha=1/2$, $\rho > 1/4$, and by arguments similar to what we did before, it leads to the Proposition \ref{highDimModel}.
\end{demo}

 \subsection{Oracle inequality: proof of Theorem 3.1}
 \label{oracle type inequality proof}
 Our oracle type inequality, Theorem 3.1, is deduced from a general model selection theorem for MLE. We state  this theorem and its proof in Section \ref{generalModelSelectionTheorem}, which is a generalization of Theorem 7.11 in \cite{MassartStFlour2007}
 for considering a random subcollection of a collection of models. Remark that this theorem could be useful in other applications, and is valid for other collections of models.
 
 Then, in Section \ref{ssPoids}, Section \ref{ssBrackets}, Section \ref{ssBernstein}, Section \ref{ssOracleInequality}, we prove the several assumptions for our specific collection of models and deduce the oracle type inequality.
\subsubsection{Model selection theorem for MLE among a random sub-collection}
\label{generalModelSelectionTheorem}
Let $(F_m)_{m \in \mathcal{M}}$ be a deterministic collection  at most countable of models.
We introduce the entropy with bracketing of the set $S$. Recall that it is defined, for every positive $\epsilon$, as the logarithm of the minimal number of brackets with 
$\dH_H$ diameter not larger than $\epsilon$ which are needed to cover $S$ and is denoted by  $\mathcal{H}_{[.]}(\epsilon,S,\dH_H)$.

In order to avoid measurability problems, we shall consider the following separability condition on the models.

\begin{assum}
\label{assumSep}
For every $m \in \mathcal{M}$, there exists some countable subset $F_m'$ of $F_m$ and a set $\mathcal{Y}$ dense in $\mathbb{R}^ p$ such that for every $f \in F_m$, there exists some sequence 
$(f_k)_{k\geq 1}$ of elements of $F'_m$ such that for every 
 $y \in \mathcal{Y}$, $\log(f_k(y))$ tends to $\log(f(y))$ as $k$ tends to infinity.
\end{assum}

\begin{theorem}
\label{ThmgeneralModelSelectionTheorem}
 Let $f^\star$ be an unknown density to be estimated from a sample of size $n$ $(\y_1,\ldots,\y_n)$. Consider $\{F_m\}_{m \in \mathcal{M}}$ some at most countable deterministic 
 collection of models where for each $m \in \mathcal{M}$, the elements of $F_m$ are assumed to be probability densities and $F_m$ fulfills Assumption \ref{assumSep}.
 Let $\{w_m\}_{m \in \mathcal{M}}$ be some family of nonnegative numbers such that:
 \begin{align}
 \label{hypK}
 \sum_{m \in \mathcal{M}} \exp({-w_m}) = \Omega < \infty . 
 \end{align}
 
We assume that for every $m \in \mathcal{M}$, 
$\sqrt{\mathcal{H}_{[.]} (\epsilon,F_m,\dH_H)}$ is integrable in $0$.
 
 Moreover, for every $m \in \mathcal{M}$, we assume that there exists $\psi_m$ on $\mathbb{R}_+$ such that $\psi_m$ is nondecreasing, $\xi \mapsto {\psi_m(\xi)}/{\xi}$ 
 is nonincreasing on $(0,+\infty)$, and
 for every $\xi \in \mathbb{R}^+$, for every $g \in F_m$, denoting by $F_m(g,\xi) = \{f\in F_m, \dH_H(f,g) \leq \xi\}$,
 \begin{align}
\label{Hm}
 \int_{0}^\xi \sqrt{\mathcal{H}_{[.]} (\epsilon,\sqrt{F_m(g,\xi)},\dH_H)}d\epsilon \leq \psi_m(\xi). 
 \end{align}

 Let $\xi_m$ be the unique positive solution of the equation $\psi_m(\xi_m)=\sqrt{n}\xi_m^2$.

 Let $\tau>0$, and for every $m \in \mathcal{M}$, let $f_m \in F_m$ such that: 
\begin{align}
\KL(f^\star,f_m) &\leq 2 \inf_{f \in F_m} \KL(f^\star,f) \nonumber ;\\
f_m &\geq \exp \left({-\tau}\right) f^\star.
\label{Bernstein}
\end{align}
 
Introduce $\{F_m\}_{m\in \tilde{\mathcal{M}}}$ some random sub-collection of $\{F_m\}_{m \in \mathcal{M}}$.

 Let $\eta \geq 0$. We consider the collection of $\eta$-maximum likelihood estimators $\{\hat{f}_{m}\}_{m \in \mathcal{\tilde{M}}}$.
 Let $pen: \mathcal{M} \rightarrow \mathbb{R}^+$, and let $\eta^{'} \geq 0$.
 We consider the $\eta'$ minimizer of the penalized criterion 
 $$\crit(m) = -\frac{1}{n} \sum_{i=1}^n \log(\hat{f}_m(\y_i))+\pen(m).$$
 Then, there exists some absolute constants $\kappa$ and $C_{\text{oracle}}$ such that wherever
 %Suppose that there exists an absolute constant $\kappa>0$ such that for all $m \in \mathcal{M}$: 
 $$\mathrm{pen}(m) \geq \kappa \left(\xi_m^2 +(1\vee \tau) \frac{w_m}{n} \right)$$
 for every $m \in \mathcal{M}$,
some random variable $\hat{m} \in \tilde{\mathcal{M}}$ such that
 $$\crit(\hat{m}) \leq \inf_{m \in \tilde{\mathcal{M}}} \crit(m) +\eta^{'}$$
exists and moreover, whatever the true density $f^\star$,
 
 $$\mathbb{E}(\dH_H^2(f^\star,\hat{f}_{\hat{m}})) \leq C_{\text{oracle}} \mathbb{E} \left( \inf_{m \in \tilde{\mathcal{M}}} \inf_{f \in F_m} \KL(f^\star,f) +\mathrm{pen}(m) \right) + (1\vee \tau)\frac{\Omega^2}{n} + \eta +\eta^{'}.$$

 \end{theorem}
 
This theorem is a generalization of Theorem $7.11$ in \cite{MassartStFlour2007} to a random subcollection of the whole collection of models.
As the proof is adapted from the proof of this theorem, we detail here only differences and we refer the interested reader to \cite{MassartStFlour2007}.
 \begin{proof}

We denote by $\gamma_n$ the empirical process and by $\bar{\gamma}_n$ the centered empirical process.
Following the proof of the Massart's theorem, easy computations lead to:
\begin{align*}
2 \KL\left(f,\frac{f+\hat{f}_{m'}}{2}\right) \leq \KL(f,f_m)+\mathrm{pen}(m)-\mathrm{pen}(m') + 2(\bar{\gamma}_n(g_m)-\bar{\gamma}_n(\hat{s}_{m'}))
\end{align*}
where
\begin{align*}
g_m = -\frac{1}{2} \log\left(\frac{f_m}{f}\right) \hspace{1cm} \text{and} \hspace{1cm} \hat{s}_m = - \log \left( \frac{f+\hat{f}_m}{2f} \right)
\end{align*}
for $m \in \tilde{\mathcal{M}}$ and $m' \in \tilde{\mathcal{M}}(m) = \left\{ m' \in \tilde{\mathcal{M}}, \gamma_n(\hat{f}_{m'}) +\mathrm{pen}(m') \leq \gamma_n(\hat{f}_m) + \mathrm{pen}(m) \right\}$.

To bound $\bar{\gamma}_n(\hat{s}_{m'})$, we use Massart's arguments.
The main difference stands in the control of $\bar{\gamma}_n(g_m)$.
As $\tilde{\mathcal{M}} \subset \mathcal{M}$ is a random subcollection of models, $\mathbb{E}(\bar{\gamma}_n(g_m)) \neq 0$.
Nevertheless, thanks to Bernstein inequality, which we may use thanks to the inequality \eqref{Bernstein}, we obtain that, for all $u>0$, with probability smaller than $\exp({-u})$,
\begin{align*}
\gamma_n(g_m) \leq \sqrt{\frac{1}{n} \alpha_\tau (1 \vee \tau) \KL(f,f_m) u} + \frac{\tau}{2n} u,
\end{align*}
where $\alpha_\tau$ is a constant depending on $\tau$.
Then, choosing $u=w_m$ for all $m \in \mathcal{M}$, where $w_m$ is defined in \eqref{hypK}, some fastidious but straightforward computations similar to those of Massart's lead to Theorem \ref{ThmgeneralModelSelectionTheorem}. 
\end{proof}

This extension has already been obtained by \cite{MaugisMeynet2012}. We remark that this is a theoretically easy extension, but quite useful in practice, {\it e.g.} for controlling large collection of models.

\subsubsection{Assumption \eqref{hypK} from Theorem \ref{ThmgeneralModelSelectionTheorem}: Construction of the weights}
\label{ssPoids}
We want to construct  a family of nonnegative numbers $\{ w_\mathbf{B}\}_{\mathbf{B} \in \mathcal{B}}$ such that
$$\sum_{\mathbf{B}\in \mathcal{B}} \exp(-w_\mathbf{B}) = \Omega < +\infty.$$
We need first to control the cardinal of $\mathcal{B}$, called the Bell number. Recall that $\mathcal{B}$ is the set of all possible partitions of the $p$ variables.
For this, we use a result of \cite{berend2010improved}, which guarantees the following inequality.
\begin{align*}
|\mathcal{B}| &\leq \left( \frac{0.792 p}{\log(p+1)} \right)^p.
\end{align*}

Then, we obtain the following result, which defines the weights needed in Assumption \eqref{hypK}.
\begin{lemma}
Let $w_{\mathbf{B}} = p \log \left(\frac{0.792 p}{\log(p+1)} \right) $.
Then, $\sum_{\mathbf{B} \in \mathcal{B}} \exp({-w_{\mathbf{B}}}) \leq 1$.
\end{lemma}

\begin{lemma}
Let $\tilde{w}_{\mathbf{B}} = p \log (p) $. Remark that $\tilde{w}_{\mathbf{B}} \geq w_{\mathbf{B}}$.
Then, $\sum_{\mathbf{B} \in \mathcal{B}} \exp({-\tilde{w}_{\mathbf{B}}}) \leq 1$.
\end{lemma}
 
\subsubsection{Assumption \eqref{Hm} from Theorem \ref{ThmgeneralModelSelectionTheorem}: Bracketing entropy}
 \label{ssBrackets}
 Let $\mathbf{B} \in \mathcal{B}$.
 Let $f \in F_{\mathbf{B}}$: it could be written $f = \Phi(0,\SSigma_{\mathbf{B}})$.
 Let $\epsilon>0$ and $\alpha>0$.
 According to Corollary \ref{constructionSBDisc2}, there exists $G \in \tilde{S}_{\mathbf{B}}^\text{disc}(\epsilon)$ such that:
 \begin{align*}
 ||\SSigma_{\mathbf{B}}-G||_2^2 \leq \frac{D_{\mathbf{B}}}{2} \wedge \epsilon.
 \end{align*}
Let $\epsilon \leq \frac{D_{\mathbf{B}}}{2}$.
For this $G$, we define the following brackets, for $\delta>0$ and $\gamma >0$, 
\begin{align*}
 u(x) &= (1+2\delta)^\gamma \phi(x|0,(1+\delta)G), \\
 l(x) &= (1+2\delta)^{-\gamma} \phi (x|0,(1+\delta)^{-1} G).
\end{align*}
According to the Proposition B.10 in \cite{Maugis2011}, if $\delta=\beta/\sqrt{3}\gamma$ and if $\epsilon = \lambda_m\beta / (3\sqrt{3}p^2)$,
the set 
$\{l,u\}$ is a $\beta$-bracket set over $F_{\mathbf{B}}$.

If we denote by $\mathcal{N}_{[.]} (\beta,F_{\mathbf{B}},\dH_H)$ the minimal number of $\epsilon$-brackets $[l,u]$ which are necessary to cover $F_{\mathbf{B}}$ and
$\mathcal{H}_{[.]} (\beta,F_{\mathbf{B}},\dH_H)$ the logarithm of this number, which corresponds to the bracketing entropy, we obtain from Corollary \ref{constructionSBDisc} and Proposition \ref{highDimModel} that
\begin{align*}
\mathcal{N}_{[.]} (\beta,F_{\mathbf{B}},\dH_H) &\leq  \left( \frac{3\sqrt{3} p^2 \lambda_M p(p-1)}{\lambda_m D_{\mathbf{B}} \beta}\right) ^{\frac{D_{\mathbf{B}}}{4}}\\
\mathcal{H}_{[.]} (\beta,F_{\mathbf{B}},\dH_H) &\leq D_{\mathbf{B}} \left( \log \mathcal{C}+\log{\left(\frac{p^3 (p-1)}{D_{\mathbf{B}}\beta}\right)}\right).
 \end{align*}
with $\mathcal{C} = \frac{3\sqrt{3}\lambda_M}{\lambda_m}$.

Subsequently, we construct $\psi_{\mathbf{B}}$ satisfying Equation \eqref{Hm}.
Here, we use the global version of the integrated square entropy rather than the local one introduced in Theorem \ref{ThmgeneralModelSelectionTheorem}.
This point is discussed in \cite{MassartStFlour2007}, Section 7.4. We should remark that no extra undesirable logarithm factor 
is obtained, since the bound is the same as the minimax lower bound.

For all $\xi >0$,
$$\int_{0}^\xi \sqrt{\mathcal{H}_{[.]} (\beta,F_{\mathbf{B}},\dH_H)} d\beta \leq \xi \sqrt{D_{\mathbf{B}} \log \mathcal{C}} + \sqrt{D_{\mathbf{B}} }\int_{0}^\xi \sqrt{\log\left(\frac{p^3 (p-1)}{D_{\mathbf{B}} \beta}\right)} d\beta.$$

According to \cite{Maugis2011},
$$\int_0^\xi \sqrt{\log\left(\frac{1}{\beta}\right)} d\beta \leq \int_0^{\xi \wedge 1} \sqrt{\log\left(\frac{1}{\beta}\right)} d\beta \leq (\xi \wedge 1) \left(\sqrt{\pi}+\sqrt{\log\left(\frac{1}{\xi \wedge 1}\right)}\right).$$

Then, denoting by $c = \sqrt{\log \mathcal{C}} + \sqrt{\pi}$, we can define $\psi_{\mathbf{B}}$ by:
$$\psi_{\mathbf{B}} (\xi) = \sqrt{D_{\mathbf{B}}} \xi \left(c+\sqrt{ \log\frac{p^3 (p-1)}{D_{\mathbf{B}}}}+\sqrt{\log{\frac{1}{\xi \wedge 1}}}\right).$$

As we want $\xi_{\mathbf{B}}$ such that $\psi_{\mathbf{B}} (\xi_{\mathbf{B}}) = \sqrt{n} \xi_{\mathbf{B}}^2$, we take:
 
\begin{align*}
\xi_{\mathbf{B}}^2 &\leq \frac{D_{\mathbf{B}}}{n} \left[2c^2 +  \log\left( \frac{p^3 (p-1)}{D_{\mathbf{B}}(\frac{D_{\mathbf{B}}}{n} c^2 \wedge 1)} \right) \right]\\
&\leq \frac{D_{\mathbf{B}}}{n} \left[2c^2 +  \log\left( \frac{p^4}{D_{\mathbf{B}}(\frac{D_{\mathbf{B}}}{n} c^2 \wedge 1)} \right) \right].
\end{align*}

 \subsubsection{Discussion on Assumption \eqref{Bernstein}}
\label{ssBernstein}
In Theorem \ref{ThmgeneralModelSelectionTheorem}, we need to do another assumption to control the randomness due to the random subcollection.
It is explained in theoretical details in the proof in Section \ref{generalModelSelectionTheorem}.
However, we would like to know if our collection of models satisfies this assumption.

Remark that the larger the parameter $\tau$, the larger the minimal penalty, and then the oracle type inequality is less accurate.

It is difficult to have a minimal convenient value of $\tau$ 
since it depends on the unknown true density $f^\star$.
Nevertheless, as $f_{\mathbf{B}}$ is expected to be close to $f^\star$, we may think that Assumption \eqref{Bernstein} is satisfied for reasonable values of $\tau$.

Remark also that if we assume that the Hellinger distance and the Kullback-Leibler divergence are equivalent, and then if we get an oracle inequality instead 
of an oracle type inequality, the Assumption \eqref{Bernstein} is satisfied.
Recall the Lemma 7.23 in \cite{MassartStFlour2007}.
\begin{lemma}
 \label{lemma::compKLdH}
 Let $P$ and $Q$ be some probability measures.
 Then
 $$\KL (P,\frac{P+Q}{2}) \geq (2 \log(2)-1) \dH_H^2(P,Q).$$
 Moreover, whenever $P<<Q$,
 \begin{align*}
  2 \dH^2_H(P,Q) \leq \KL(P,Q) \leq 2\left(2+\log\left(\left|\left|\frac{dP}{dQ}\right|\right|_\infty\right)\right) \dH^2_H(P,Q).
 \end{align*}

\end{lemma}

Another point of view could be used: we have already assumed that the true distribution is absolutely continuous with respect to the Lebesgue measure, and that there exists a density $f^ \star$.
 To satisfy Assumption \eqref{Bernstein}, we may assume that the densities have the same support $\mathbb{R}$, and that the covariance matrix of $f^ \star$ is bounded.
In this case, the constant $\tau$ will depend on those bounds.

\subsubsection{Oracle type inequality}
\label{ssOracleInequality}

According to Theorem \ref{ThmgeneralModelSelectionTheorem}, there exists absolute constants $\kappa$ and $C_{\text{oracle}}$ such that wherever 
$$\pen(\mathbf{B}) \geq \kappa \frac{D_\mathbf{B}}{n}\left(2c^2 + \log\left( \frac{p^4}{D_\mathbf{B} (\frac{D_\mathbf{B}}{n} c^2 \wedge 1)}\right) + \frac{1 \vee \tau}{n} p \log\left( \frac{0.792 p}{\log(p+1)}\right) \right)$$
for every $\mathbf{B} \in \mathcal{B}$, some random variable $\hat{\mathbf{B}} \in \mathcal{B}_\Lambda$ such that 
$$\hat{\mathbf{B}} = \underset{\mathbf{B} \in \mathcal{B}_\Lambda}{\operatorname{argmin}} \crit(\mathbf{B})$$
exists and, moreover, whatever the true density $f^\star$,
$$\mathbb{E}(\dH_H^2(f^\star,\hat{f}_{\hat{\mathbf{B}}})) \leq C_{\text{oracle}} (\inf_{\mathbf{B} \in \mathcal{B_\Lambda}} \inf_{f \in F_\mathbf{B}}\KL(f^\star,f) +\pen(\mathbf{B})) +\frac{1\vee \tau}{n}.$$

In these graphical models, we remark that every minimizer is tractable (the maximum likelihood estimator and the minimizer of the criterion $\crit$), then we let $\eta=\eta'=0$.

Moreover, the weights do not depend on the model $\mathbf{B}$, then we can rewrite the oracle type inequality.

There exists absolute constants $\kappa$ and $C_{\text{oracle}}$ such that whenever 
$$\pen(\mathbf{B}) \geq \kappa \frac{D_\mathbf{B}}{n}\left(2c^2 + \log\left( \frac{p^4}{D_\mathbf{B} (\frac{D_\mathbf{B}}{n} c^2 \wedge 1)}\right) \right)$$
for every $\mathbf{B} \in \mathcal{B}$, some random variable $\hat{\mathbf{B}} \in \mathcal{B}_\Lambda$ such that 
$$\hat{\mathbf{B}} = \underset{\mathbf{B} \in \mathcal{B}_\Lambda}{\operatorname{argmin}} \crit(\mathbf{B})$$
exists and, moreover, whatever the true density $f^\star$,
$$\mathbb{E}(\dH_H^2(f^\star,\hat{f}_{\hat{\mathbf{B}}})) \leq C_{\text{oracle}} (\inf_{\mathbf{B} \in \mathcal{B}_\Lambda} \inf_{f \in F_\mathbf{B}}\KL(f^\star,f) +\pen(\mathbf{B})) +
\frac{1\vee \tau}{n} p \log(p).$$

This version is slightly better, because the penalty does not depend on $\tau$ which is related to the unknown true density.
In that case, the penalty is  explicitly defined through some tractable constants depending only on the definition of the collection of models.

\subsection{Lower bound for the minimax risk: proof of Theorem 3.2}
\label{lower bound proof}

Fix $\mathbf{B} \in \mathcal{B}$.

%We will show that 
% 
% \begin{align}
% \label{inegMinimax1}
%\max_{f \in F_{\mathbf{B}}(r)} \mathbb{E} (\dH_H^2(\hat{f}_{\mathbf{B}},f)) &\geq C \frac{D_{\mathbf{B}}}{n} \left(r^2 \wedge \left(1+\log \frac{p(p-1)}{2D_{\mathbf{B}}} \right)\right).
% \end{align}

%Let $F_{\mathbf{B}}^\text{disc}$ be the discrete space constructed in Corollary \ref{constructionFBDisc}.
Let $f^\star=\phi(0,{\SSigma}^\star)$ be the true density.
Let $\hat{f}$ be the considered estimator.
Let $\mathcal{F} \subset F_\mathbf{B}$ a subset of densities.
We define $\tilde{f} = \underset{f \in \mathcal{F}}{\operatorname{argmin}} \left\{ \dH_H(\hat{f},f) \right\}$.

First, we have:
\begin{align}
\label{inegTriang}
\dH_H(f^\star,\tilde{f}) \leq \dH_H(f^\star,\hat{f}) + \dH_H(\hat{f},\tilde{f}) \leq 2 \dH_H(\hat{f},f^\star). 
\end{align}

Secondly, we have:
\begin{align}
\dH_H(\tilde{f},f^\star)^2 &\geq \mathbf{1}_{f^\star \neq \tilde{f}} \min_{f \neq f^\star} \dH_H(f^\star,f)^2 \nonumber\\
\mathbb{E}(\dH_H(\tilde{f},f^\star)^2) &\geq \mathbb{P}_{f^\star}(f^\star \neq \tilde{f}) \min_{f\neq f^\star} \dH_H(f^\star,f)^2. 
\label{proba}
\end{align}

Then, by combining \eqref{inegTriang} and \eqref{proba} we obtain:
\begin{align}
\label{eq0}
\max_{f^\star \in \mathcal{F}} \mathbb{E}(\dH_H^2(\hat{f},f^\star)) \geq \frac{1}{4} \max_{f^\star \in \mathcal{F}} \left[ \mathbb{P}_f^\star (f^\star \neq \tilde{f}) \min_{f \neq f^\star} 
\dH_H^2(f^\star,f) \right].
\end{align}

We need to design a lower bound for $\max_{f \in \mathcal{F}} \mathbb{P}_f (f \neq \tilde{f}).$
For this purpose, we use the Birg\'e's lemma.

\begin{lemma}[Birg\'e's lemma]
\label{BirgeLemma}
 Let $(\mathbb{P}_f)_{f \in \mathcal{F}}$ a probability family, and $(A_f)_{f\in \mathcal{F}}$ some event pairwise disjoints.
 Then,
 $$\min_{f \in \mathcal{F}} \mathbb{P}_f(A_f) \leq \frac{2e}{2e+1} \vee \frac{\max_{f\in \mathcal{F}} \KL(\mathbb{P}_f,\mathbb{P}_0)}{\log(1+ \text{card} (\mathcal{F}))}.$$
 
\end{lemma}

If we use Birg\'e's lemma to control $\max_{f \in F_{\mathbf{\mathbf{B}}}(r)} \mathbb{P}_f (f \neq \tilde{f})$ in \eqref{eq0}, we obtain:
\begin{align}
\label{eq01}
\max_{f^\star \in \mathcal{F}} \mathbb{E}(\dH_H^2(\hat{f},f^\star)) \geq \frac{1}{4(2e+1)} \max_{f^\star \in \mathcal{F}} (\min_{f \neq f^\star} \dH_H^2(f^\star,f))
\end{align}
if
$$ \frac{2e}{2e+1} \geq \frac{\max_{f\in \mathcal{F}} \KL(\mathbb{P}_f,\mathbb{P}_0)}{\log(1+ \text{card} (\mathcal{F}))}$$

Then, we distinguish between low and high dimensional models, using the discretization constructed in Section \ref{Discretization} instead of $\mathcal{F}$.

\textbf{First case: $p(p-1)/2 \geq 4D_{\mathbf{B}}$}

Let $\epsilon = D_{\mathbf{B}}/2$. Let $\tilde{S}_{\mathbf{B}}^{\text{disc}}(D_{\mathbf{B}}/2)$ the discrete space constructed in Corollary \ref{constructionSBDisc2}, and the following quantity for $r>0$:
\begin{align*}
 \tilde{F}_{\mathbf{B}}(r) &= \left\{ r S, S \in \tilde{S}_{\mathbf{B}}^{\text{disc}}\left(\frac{D_{\mathbf{B}}}{2}\right) \right\}.
\end{align*}
Using Lemma \ref{lemma::compKLdH} for comparing the $L_2$ norm and the Hellinger distance, we get
\begin{align}
\label{eq1}
\max_{f^\star \in \tilde{F}_{\mathbf{B}}(r)} \mathbb{E}(\dH_H^2(\hat{f},f^\star)) \geq \frac{1}{4(2e+1)} \frac{1}{4+p/2 \log(\lambda_M/\lambda_m)} \frac{1}{2} \frac{\lambda_m p^3}{e_m^2} D_{\mathbf{B}} r^2
\end{align}

if the following inequality is satisfied:
\begin{align}
\label{eq2}
 \max_{f_1,f_2 \in \tilde{F}_{\mathbf{B}}(r)} (n \KL(f_1,f_2)) \leq \frac{2e}{2e+1} \log\left(1+\text{card}(\tilde{F}_{\mathbf{B}}(r))\right).
\end{align}

The inequality \eqref{eq2} is satisfied if the inequality \eqref{eq3} is fulfilled, with:

\begin{align}
\label{eq3}
 \frac{n}{2} p^3 \frac{\lambda_m}{e_m^2} D_{\mathbf{B}} r^2 \leq \frac{2e}{2e+1} \rho D_{\mathbf{B}} \log \left(\frac{p(p-1)2 \lambda_M}{ D_{\mathbf{B}}^2} \right) .
\end{align}

Then, we can replace this condition in \eqref{eq1} and we obtain:

 \begin{align*}
\max_{f^\star \in \tilde{F}_{\mathbf{B}}(r)} \mathbb{E} (\dH_H^2(\hat{f}_{\mathbf{B}},f^\star)) &\geq \tilde{C}_1 \frac{D_{\mathbf{B}}}{n}\log \frac{\tilde{C}_2}{D_{\mathbf{B}}^2}
 \end{align*}

 with
 $$\tilde{C}_1= \frac{2e}{4(2e+1)^2} \frac{1}{4+p/2 \log(\lambda_M/\lambda_m)} \rho,$$
 $\tilde{C}_2 = 2p(p-1)\lambda_M$
 and  $ 0.233 \leq \rho \leq 0.234$.
 
\textbf{Second case: $p(p-1)/2 \leq 4D_{\mathbf{B}}$}
Let $\epsilon = p(p-1)/16$.
We use the Proposition \ref{highDimModel}, and consider $\check{S}_\mathbf{B}^{\text{disc}}(p(p-1)/16)$ and construct $\check{F}_{\mathbf{B}}(r)$ as previously.
Then, we get the following.

\begin{align*}
 \max_{f^\star \in \tilde{F}_{\mathbf{B}}(r)} \mathbb{E}(\dH_H^2(\hat{f},f^\star) \geq \frac{1}{4(2e+1)}\frac{1}{4+p/2 \log(\lambda_M/\lambda_m)} \frac{1}{2} \frac{\lambda_m p^3}{e_m^2} D_{\mathbf{B}} r^2,
\end{align*}
if
\begin{align*}
 \frac{n}{2}  \frac{\lambda_m p^3}{e_m^2} D_{\mathbf{B}} r^2 \leq \frac{2e}{2e+1}  \frac{D_{\mathbf{B}}}{4} \log\left(\frac{2 \lambda_M}{p(p-1)/16}\right).
\end{align*}

Then, we obtain the following bound:
\begin{align*}
 \max_{f^\star \in \tilde{F}_{\mathbf{B}}(r)} \mathbb{E}(\dH_H^2(\hat{f},f^\star) \geq \frac{1}{4(2e+1)}\frac{1}{4+p/2 \log(\lambda_M/\lambda_m)}  \frac{2e}{2e+1}  \frac{D_{\mathbf{B}}}{4n} \log\left(\frac{2 \lambda_M}{p(p-1)/16}\right).
\end{align*}

\textbf{Conclusion}
As $\tilde{F}_{\mathbf{B}}(r) \subset F_{\mathbf{B}}$, and $\check{F}_{\mathbf{B}}(r) \subset F_{\mathbf{B}}$, choosing $r=(1+\log(2 p(p-1)\lambda_M/D_{\mathbf{B}}^2))^{1/2}$, we get that: 
 \begin{align*}
\max_{f^\star \in F_{\mathbf{B}}} \mathbb{E} (\dH_H^2(\hat{f}_{\mathbf{B}},f^\star)) &\geq C_{\text{minim}} \frac{D_{\mathbf{B}}}{n}\left(1+\log \frac{2 p(p-1)\lambda_M}{D_{\mathbf{B}}^2} \right),
 \end{align*}

 with: 
 $$C_{\text{minim}}= \frac{e}{(2e+1)^2} \frac{1}{8+p\log(\lambda_M/\lambda_m)} \rho\left( \frac{1}{4} \log \left(\frac{32 \lambda_M}{p(p-1)}\right) \wedge 1 \right),$$
 and with $ 0.233 \leq \rho \leq 0.234$.

 \bibliographystyle{alpha}

\bibliography{draft_EDMG_arxiv}

\newcommand{\etalchar}[1]{$^{#1}$}
\begin{thebibliography}{RMRMMC15}

\bibitem[ACM09]{2009_EJS_Chiquet}
C.~Ambroise, J.~Chiquet, and C.~Matias.
\newblock Inferring sparse {Gaussian} graphical models with latent structure.
\newblock {\em Electronic Journal of Statistics}, 3:205--238, 2009.

\bibitem[AL13]{Allen2013}
G.~Allen and Z.~Liu.
\newblock {A local Poisson graphical model for inferring networks from
  sequencing data.}
\newblock {\em {IEEE Transactions on NanoBioscience}}, 12(3):189--198, 2013.

\bibitem[AM09]{Arlot2010}
S.~Arlot and P.~Massart.
\newblock Data-driven calibration of penalties for least-squares regression.
\newblock {\em Journal of Machine Learning Research}, 10:245--279, 2009.

\bibitem[ANW{\etalchar{+}}14]{Akbani2014pan}
R.~Akbani, P.~K.~S. Ng, H.~MJ Werner, M.~Shahmoradgoli, F.~Zhang, Z.~Ju,
  W.~Liu, J-Y Yang, K.~Yoshihara, J.~Li, et~al.
\newblock A pan-cancer proteomic perspective on the cancer genome atlas.
\newblock {\em Nature communications}, 5, 2014.

\bibitem[BCJ15]{Bouveyron2015}
C.~Bouveyron, E.~C\^ome, and J.~Jacques.
\newblock The discriminative functional mixture model for a comparative
  analysis of bike sharing systems.
\newblock {\em The Annals of Applied Statistics}, in press, 2015.

\bibitem[BEGd08]{Banerjee2008}
O.~Banerjee, L.~El~Ghaoui, and A.~d'Aspremont.
\newblock Model selection through sparse maximum likelihood estimation for
  multivariate gaussian or binary data.
\newblock {\em Journal of Machine Learning Research}, 9:485--516, 2008.

\bibitem[BGH09]{Baraud2009}
Y.~Baraud, C.~Giraud, and S.~Huet.
\newblock Gaussian model selection with an unknown variance.
\newblock {\em The Annals of Statistics}, 37(2):630--672, 2009.

\bibitem[Bir05]{Birge2005}
L.~Birg\'e.
\newblock A new lower bound for multiple hypothesis testing.
\newblock {\em Information Theory, IEEE Transactions}, 51(4):1611--1615, 2005.

\bibitem[BL08]{Bickel2008}
P.J. Bickel and E.~Levina.
\newblock Regularized estimation of large covariance matrices.
\newblock {\em The Annals of Statistics}, 36(1):199--227, 2008.

\bibitem[BM01]{BirgeMassart2001}
L.~Birgé and P.~Massart.
\newblock Gaussian model selection.
\newblock {\em Journal of the European Mathematical Society}, 3(3):203--268,
  2001.

\bibitem[BM07]{BirgeMassart2007}
L.~Birg{\'{e}} and P.~Massart.
\newblock {Minimal penalties for {G}aussian model selection}.
\newblock {\em Probability Theory \& Related Fields}, 138(1-2), 2007.

\bibitem[BMM12]{Baudry2012}
J-P Baudry, C.~Maugis, and B.~Michel.
\newblock Slope heuristics: overview and implementation.
\newblock {\em Statistics and Computing}, 22(2):455--470, 2012.

\bibitem[BT10]{berend2010improved}
D.~Berend and T.~Tassa.
\newblock Improved bounds on bell numbers and on moments of sums of random
  variables.
\newblock {\em Probability and Mathematical Statistics}, 30(2):185--205, 2010.

\bibitem[CN06]{igraph}
Gabor Csardi and Tamas Nepusz.
\newblock The igraph software package for complex network research.
\newblock {\em InterJournal}, Complex Systems:1695, 2006.

\bibitem[CZZ10]{Cai2010}
T.~Cai, C-H Zhang, and H.~Zhou.
\newblock Optimal rates of convergence for covariance matrix estimation.
\newblock {\em The Annals of Statistics}, 38(4):2118--2144, 2010.

\bibitem[DWW14]{Danaher2014}
P.~Danaher, P.~Wang, and D.M. Witten.
\newblock The joint graphical lasso for inverse covariance estimation across
  multiple classes.
\newblock {\em Journal of the Royal Statistical Society: Series B (Statistical
  Methodology)}, 76(2):373--397, 2014.

\bibitem[FHT08]{Friedman2007}
J.~Friedman, T.~Hastie, and R.~Tibshirani.
\newblock Sparse inverse covariance estimation with the graphical lasso.
\newblock {\em Biostatistics}, 9(3):432--441, 2008.

\bibitem[FLL11]{Frazee2011}
A.~C. Frazee, B.~Langmead, and J.~T. Leek.
\newblock {ReCount: a multi-experiment resource of analysis-ready RNA-seq gene
  count datasets}.
\newblock {\em BMC Bioinformatics}, 12(449), 2011.

\bibitem[GB09]{mvtnorm}
A.~Genz and F.~Bretz.
\newblock {\em Computation of Multivariate Normal and T Probabilities}.
\newblock Springer Publishing Company, Incorporated, 1st edition, 2009.

\bibitem[GHV12]{GiraudHuetVerzelen2012}
C.~Giraud, Sylvie Huet, and Nicolas Verzelen.
\newblock {Graph selection with GGMselect}.
\newblock {\em Statistical Applications in Genetics and Molecular Biology},
  11(3), 2012.

\bibitem[Gir08]{Giraud2008}
C.~Giraud.
\newblock {Estimation of Gaussian graphs by model selection}.
\newblock {\em Electronic Journal of Statistics}, 2:542--563, 2008.

\bibitem[GLMZ11]{Guo2011}
J.~Guo, E.~Levina, G.~Michailidis, and J.~Zhu.
\newblock {Joint estimation of multiple graphical models.}
\newblock {\em Biometrika}, 98(1):1--15, 2011.

\bibitem[GW00]{Genovese2000}
C.~Genovese and L.~Wasserman.
\newblock Rates of convergence for the {G}aussian mixture sieve.
\newblock {\em The Annals of Statistics}, 28(4):1105--1127, 2000.

\bibitem[HA85]{Hubert1985}
L.~Hubert and P.~Arabie.
\newblock Comparing partitions.
\newblock {\em Journal of Classification}, 2(1):193--218, 1985.

\bibitem[HSDR14]{Hsieh2014}
Cho-Jui Hsieh, M\'{a}ty\'{a}s~A. Sustik, Inderjit~S. Dhillon, and Pradeep
  Ravikumar.
\newblock Quic: Quadratic approximation for sparse inverse covariance
  estimation.
\newblock {\em Journal of Machine Learning Research}, 15:2911--2947, 2014.

\bibitem[HSNP15]{Hyodo2015}
M.~Hyodo, N.~Shutoh, T.~Nishiyama, and T.~Pavlenko.
\newblock Testing block-diagonal covariance structure for high-dimensional
  data.
\newblock {\em Statistica Neerlandica}, 69(4):460--482, 2015.

\bibitem[KSI{\etalchar{+}}11]{Krumsiek2011}
J.~Krumsiek, K.~Suhre, T.~Illig, J.~Adamski, and F.~J Theis.
\newblock {Gaussian graphical modeling reconstructs pathway reactions from
  high-throughput metabolomics data}.
\newblock {\em BMC Systems Biology}, 5(1):21, 2011.

\bibitem[Leb05]{Lebarbier2005}
E.~Lebarbier.
\newblock Detecting multiple change-points in the mean of gaussian process by
  model selection.
\newblock {\em Signal Processing}, 85(4):717 -- 736, 2005.

\bibitem[Mas07]{MassartStFlour2007}
P.~Massart.
\newblock {\em Concentration inequalities and model selection}.
\newblock Lecture Notes in Mathematics. Springer, 33, 2003, Saint-Flour,
  Cantal, 2007.

\bibitem[MB06]{MeinshausenBuhlmann}
N.~Meinshausen and P.~B{\"u}hlmann.
\newblock High-dimensional graphs and variable selection with the lasso.
\newblock {\em The Annals of Statistics}, 34(3):1436--1462, 2006.

\bibitem[MH12]{Mazumder2012}
R.~Mazumder and T.~Hastie.
\newblock {Exact covariance thresholding into connected components for
  large-scale Graphical Lasso}.
\newblock {\em Journal of Machine Learning Research}, 13:781--794, 2012.

\bibitem[MM11]{Maugis2011}
C.~Maugis and B.~Michel.
\newblock A non asymptotic penalized criterion for {G}aussian mixture model
  selection.
\newblock {\em ESAIM. Probability and Statistics.}, 15:41--68, 2011.

\bibitem[MMR12]{MaugisMeynet2012}
C.~Meynet and C.~Maugis-Rabusseau.
\newblock {A sparse variable selection procedure in model-based clustering}.
\newblock Research report, Department of Mathematics, Université Paris-Sud,
  hal-00734316, 2012.

\bibitem[PBT12]{Pavlenko2012}
T.~Pavlenko, A.~Bj{\"o}rkstr{\"o}m, and A.~Tillander.
\newblock Covariance structure approximation via glasso in high dimensional
  supervised classification.
\newblock {\em Journal of Applied Statistics}, 39(8):1643--1666, 2012.

\bibitem[Pic10]{Pickrell2010}
J~Pickrell.
\newblock {Understanding mechanisms underlying human gene expression variation
  with RNA sequencing}.
\newblock {\em Nature}, 464(7289):768--772, 2010.

\bibitem[RGCJ13]{Rau2013}
A.~Rau, M.~Gallopin, G.~Celeux, and F.~Jaffr\'{e}zic.
\newblock {Data-based filtering for replicated high-throughput transcriptome
  sequencing experiments}.
\newblock {\em Bioinformatics}, 29(17):2146--2152, 2013.

\bibitem[RMRMMC15]{Rau2015}
A.~Rau, C.~Maugis-Rabusseau, M-L Martin-Magniette, and G.~Celeux.
\newblock {Co-expression analysis of high-throughput transcriptome sequencing
  data with Poisson mixture models}.
\newblock {\em Bioinformatics}, 31(9):1420--1427, 2015.

\bibitem[TWS15]{Tan2015}
K.~Tan, D.M. Witten, and A.~Shojaie.
\newblock {The Cluster Graphical Lasso for improved estimation of Gaussian
  graphical models}.
\newblock {\em Computational Statistics \& Data Analysis}, 85:23--36, 2015.

\bibitem[Ver12]{Verzelen2012}
N.~Verzelen.
\newblock {Minimax risks for sparse regressions: Ultra-high dimensional
  phenomenons}.
\newblock {\em Electronic Journal of Statistics}, 6:38--90, 2012.

\bibitem[WFS11]{Witten2011GM}
D.~M. Witten, J.~H. Friedman, and N.~Simon.
\newblock {New insights and faster computations for the Graphical Lasso}.
\newblock {\em Journal of Computational and Graphical Statistics},
  20(4):892--900, 2011.

\bibitem[Whi90]{Whittaker1990}
J.~Whittaker.
\newblock {\em Graphical Models in Applied Multivariate Statistics}.
\newblock Wiley Publishing, 1990.

\bibitem[YL07]{Yuan2007}
M.~Yuan and Y.~Lin.
\newblock {Model selection and estimation in the Gaussian graphical model}.
\newblock {\em Biometrika}, 94(1):19--35, 2007.

\bibitem[YL11]{Yin2011}
J.~Yin and H.~Li.
\newblock A sparse conditional gaussian graphical model for analysis of
  genetical genomics data.
\newblock {\em The Annals of Applied Statistics}, 5(4):2630--2650, 2011.

\bibitem[ZLR{\etalchar{+}}12]{Zhao2012}
T.~Zhao, H.~Liu, K.~Roeder, J.~Lafferty, and L.~Wasserman.
\newblock The huge package for high-dimensional undirected graph estimation in
  {R}.
\newblock {\em The Journal of Machine Learning Research}, 13(1):1059--1062,
  2012.

\end{thebibliography}

%\section{First section}
%Here is some text for the first section, and a label\label{sec1}.
%Uses version \RRfileversion\ of the package.\newpage
%\section{Second section}
%Text for the second section. This is closely related to the text in
%section \ref{sec1} on page \pageref{sec1}. \newpage
%\tableofcontents

\end{document}